\documentclass[11pt,reqno]{amsart}
\usepackage{mathtools}
\usepackage{amsmath,amssymb,amsthm}
\usepackage{fullpage}
\usepackage{mathrsfs,bm,color}
\usepackage{hyperref}
\usepackage{bbm}
\usepackage{cite}
\usepackage{xcolor}

\hypersetup{
	colorlinks,
	linkcolor={red!50!black},
	citecolor={green!50!black},
	urlcolor={blue!80!black}
}

\theoremstyle{definition}

\theoremstyle{plain}
\newtheorem{theo}{Theorem}
\newtheorem*{theo*}{Theorem}
\newtheorem{cor}{Corollary}

\newtheorem{lem}{Lemma}

\numberwithin{theo}{section}
\numberwithin{lem}{section}
\numberwithin{cor}{section}
\numberwithin{proposition}{section}
\numberwithin{equation}{section}
\numberwithin{figure}{section}

\noeqref{SHOWLABEL}

\newcommand{\mmm}{\mathrel{}\mid\mathrel{}}

\newcommand{\bs}{\bm{\sigma}}
\newcommand{\iii}{i_1,\dots,i_p}
\newcommand{\E}{\mathbb{E}}

\newcommand{\sO}{\mathcal{O}}

\newcommand{\cL}{\mathcal{L}}

\newcommand{\sP}{\mathscr{P}}

\newcommand{\N}{\mathbb{N}}

\newcommand{\R}{\mathbb{R}}

\newcommand{\tr}{\operatorname{tr}}
\newcommand{\iid}{\text{i.i.d. }}

\DeclareMathOperator{\pP}{\mathbb{P}}

\DeclareMathOperator{\1}{\mathbbm{1}}
\DeclareMathOperator{\Cov}{\operatorname{Cov}}
\newcommand{\eqdist}{\stackrel{d}{=}}
\newcommand{\bt}{\bm{\tau}}
\newcommand{\bo}{\bm{\omega}}

\newcommand{\bvs}{\bm{\sigma}}
\newcommand{\bvp}{\bm{\rho}}
\newcommand{\bvz}{\vec{z}}
\newcommand{\bvt}{\bm{\tau}}

\newcommand{\bto}{\tilde{\bm{\omega}}}
\newcommand{\bvr}{\bm{\rho}}
\newcommand{\bvo}{\bm{\omega}}
\newcommand{\vs}{\sigma}
\newcommand{\vb}{\beta}
\newcommand{\vh}{h}
\newcommand{\veh}{\vec{h}}
\newcommand{\vo}{\omega}

\newcommand{\vnu}{\nu}
\newcommand{\vx}{x}
\newcommand{\vp}{\rho}
\newcommand{\vg}{g}
\newcommand{\vt}{\tau}

\newcommand{\bvb}{\bm{\beta}}
\newcommand{\ubvs}{\tilde{\bm{\sigma}}}
\newcommand{\ubvp}{\tilde{\bm{\rho}}}

\newcommand{\ubvo}{\tilde{\bm{\omega}}}
\newcommand{\bhs}{\hat{\bm{\sigma}}}

\newcommand{\bR}{\bm{R}}
\newcommand{\btR}{\hat{\bm{R}}}
\newcommand{\bQ}{\bm{Q}}
\newcommand{\bA}{\bm{A}}
\newcommand{\bL}{\bm{\Lambda}}
\newcommand{\bD}{\bm{\Delta}}
\newcommand{\bd}{\bm{D}}
\newcommand{\bP}{\bm{U}}
\newcommand{\bC}{\bm{C}}
\newcommand{\bzero}{\bm{0}}

\newcommand{\bB}{\bm{B}}

\newcommand{\bI}{\bm{I}}
\newcommand{\bxi}{\bm{\xi}}
\newcommand{\btheta}{\bm{\theta}}
\newcommand{\Sum}{\mathrm{Sum}}

\newcommand{\trans}{\mathsf{T}}

\renewcommand{\phi}{\varphi}
\renewcommand{\epsilon}{\varepsilon}

\newenvironment{linsys}[2][m]{%
	\setlength{\arraycolsep}{.1111em} 
	\begin{array}[#1]{@{}*{#2}{rc}r@{}} 
	}{%
\end{array}}

\title[Free Energy of Multiple Systems of Spherical Spin Glasses with Constrained Overlaps]
{Free Energy of Multiple Systems of Spherical Spin Glasses with Constrained Overlaps}

\author{Justin Ko}
\address[Justin Ko]{Department of Mathematics, University of Toronto, Partially Supported by NSERC}
\email{jko@math.toronto.edu}

\date{\today}

\begin{document}

\begin{abstract}
The free energy of multiple systems of spherical spin glasses with constrained overlaps was first studied in \cite{PTSPHERE}. The authors proved an upper bound of the constrained free energy using Guerra's interpolation. In this paper, we prove this upper bound is sharp. Our approach combines the ideas of the Aizenman--Sims--Starr scheme in \cite{CASS} and the synchronization mechanism used in the vector spin models in \cite{PVS} and \cite{PPotts}. We derive a vector version of the Aizenman--Sims--Starr scheme for spherical spin glass and use the synchronization property of arrays obeying the overlap-matrix form of the Ghirlanda--Guerra identities to prove the matching lower bound. 
\end{abstract}

\maketitle

\section{Introduction}

In \cite{TSPHERE}, Talagrand proved a formula for the free energy of the spherical mixed even-$p$-spin model originally considered by Crisanti and Sommers in \cite{crisanti1992sphericalp}. It was later extended to general mixed $p$-spin models by Chen in \cite{CASS}. This formula is the analogue of the classical Parisi formula for the Sherrington--Kirkpatrick model \cite{parisi1979infinite,parisi1980sequence} proved in \cite{talagrand2006parisi}. 

This paper is on the free energy of multiple copies of spherical spin glasses with constrained overlaps. The free energy of this model was studied in \cite{PTSPHERE}, where an analogue of the Guerra replica symmetry breaking bound \cite{guerra2003broken} was derived and used in several applications. The goal of this paper is to prove that this upper bound is sharp.

There are several motivations for this paper. In \cite{arous2017spectral}, spectral gap estimates for generic spherical models were proved under various conditions on the Parisi measure. Our free energy formulas can be used to possibly prove large deviation principles to extend these spectral gap estimates to the larger class of mixed even-$p$-spin spherical models. Another application of the free energy formula is possibly proving that chaos in temperature in some full-RSB spherical models cannot be detected at the level of the free energy, as was predicted in \cite{rizzo2002ultrametricity} and recently proven geometrically in \cite{subag2018free}. See \cite{chen2017temperature} for some related results on temperature chaos for spherical models.

The main tool that allows us to prove the matching lower bound is the overlap synchronization mechanism developed by Panchenko in \cite{panchenko2015free,PVS,PPotts} to study multi-species models and models with vectors spins. This mechanism is a consequence of the ultrametric structure of generalized overlaps that satisfy the Ghirlanda--Guerra identities \cite{GG,guerra1996overlap} which was proved in \cite{PUltra}. Synchronization was used recently in other contexts in \cite{aukoshsubhabratajustin2017,contucci2018multi}, and in this paper we give another application. Besides this, our proof is based on a variant of the Aizenman--Sims--Starr scheme for spherical models developed in \cite{CASS}.

Lastly, we refer the reader to \cite{jagannath2018bounds,jagannath2017low,chen2015disorder,auffinger2017energy,gheissari2016spectral,subag2017geometry,chen2017parisi} for other recent work where various aspects of the spherical models have been studied.

\section{Model Description}
Fix $n \geq 1$. The main goal is to find a formula for the free energy of $n$ constrained copies of spherical spin glasses. The copies are coupled by constraining their overlaps and can possibly exist at different temperatures. We start by introducing the usual spherical spin glass model.

Let $S_N$ be the sphere in $\R^N$ of radius $\sqrt{N}$ and denote the configuration of the $j$th copy by
\begin{equation}\label{SHOWLABEL}
\bs(j) = \big(\vs_1(j), \dots, \vs_N(j)\big) \in S_N.
\end{equation}
For $p\geq 2$, the $p$-spin Hamiltonian is denoted by
\begin{equation}\label{SHOWLABEL}
H_{N,p}(\bs(j) ) = \frac{1}{N^{(p-1)/2}} \sum_{1 \leq \iii \leq N} g_{\iii} \vs_{i_1}(j) \cdots \vs_{i_p}(j),
\end{equation}
where $g_{\iii}$ are \iid standard Gaussian for all $p\geq 2$ and indices $(\iii)$. The corresponding even mixed $p$-spin Hamiltonian for the $j$th copy at temperature $(\vb_p(j))_{p \geq 2}$ is denoted by
\begin{equation}\label{eq:hammixp}
H^j_N(\bvs) = \sum_{p \geq 2} \vb_p(j) H_{N,p}(\bs(j) ).
\end{equation}
We assume that the inverse temperatures satisfy $\sum_{p \geq 2} 2^p \vb_p^{\, 2}(j) < \infty$ for all $j \leq n$, so that \eqref{eq:hammixp} is well-defined, and that $\vb_p(j) = 0$ for odd $p$. 

We now introduce the model for a system of $n$ copies of spherical spin glass. A configuration of $n$ copies can be viewed as vector spins,
\begin{equation}\label{SHOWLABEL}
\bvs = (\vs_1, \dots, \vs_N) \in (\R^n)^N,
\end{equation}
where the vector entries of $\bvs$ are denoted by
\begin{equation}\label{SHOWLABEL}
\vs_i = \big( \vs_i(1), \dots, \vs_i(n) \big) \in \R^n.
\end{equation}
The configurations $\bvs$ are restricted to the set
\begin{equation}\label{SHOWLABEL}
S_N^n = \big\{ \bvs \in (\R^N)^n \mmm \| \bs(j)\| = \sqrt{N} \text{ for all } j \leq n \big\},
\end{equation}
where $\| \cdot \|$ is the Euclidean norm on $\R^N$. The Hamiltonian of $n$ copies of even mixed $p$-spin models of spherical spin glasses is denoted by
\begin{equation}\label{eq:hamiltonian}
H_N(\bvs) = \sum_{j \leq n} H^j_N(\bvs).
\end{equation}

The upper indices $\ell \geq 1$ of the configurations $\bvs^\ell$ index sequences of spin configurations. The Hamiltonian is a Gaussian process indexed by $\bvs^\ell \in S_N^n$ with covariance given by functions of normalized inner products. The inner products, or overlaps, of the configurations of copy $\bs^{\ell}(j)$ and $\bs^{\ell'}(j')$ is denoted by 
\begin{equation}\label{SHOWLABEL}
R_{\ell,\ell'}^{j,j'} = R_{\ell,\ell'}^{j,j'}\bigl( \bvs^\ell(j), \bvs^{\ell'}(j')\bigr) = \frac{1}{N} \sum_{i \leq N} \vs^\ell_i(j) \vs^{\ell'}_i(j').
\end{equation}
The overlaps of vector configurations $\bvs^\ell$ and $\bvs^{\ell'}$ are given by the overlap matrices
\begin{equation}\label{SHOWLABEL}
\bm{R}_{\ell, \ell'} = \bR(\bvs^\ell, \bvs^{\ell'}) = \big(R_{\ell,\ell'}^{j,j'}\big)_{j,j' \leq n} = \frac{1}{N}\sum_{i \leq N} \vs_i^{\ell} \otimes \vs_i^{\ell'}.
\end{equation}
The overlaps are always normalized by the dimension of the vectors in the inner product. Let $x \in \R^n$ and let $\bA = (A_{j,j'})_{j,j' \leq n} \in \R^{n \times n}$. Consider the real valued convex function
\begin{equation}\label{SHOWLABEL}
\xi_{j,j'}(x) = \sum_{p \geq 2} \vb_p(j)\vb_p(j') x^p
\end{equation}
and its matrix valued counterpart
\begin{equation}\label{eq:ximatrix}
\bxi(\bA) = \big(\xi_{j,j'}(A_{j,j'})\big)_{j,j' \leq n} = \sum_{p \geq 2} (\vb_p \otimes \vb_p) \odot \bA^{\circ p},
\end{equation}
where $\otimes$ is the outer product on vectors in $\R^n$ and $\odot$ is the Hadamard product on $n \times n$ matrices. It is easy to check that the mixed $p$-spin Hamiltonian of the copies \eqref{eq:hammixp} are centered Gaussian processes with covariance
\begin{equation}\label{SHOWLABEL}
\E H^j_N\bigl(\bvs^{\ell}\bigr) H^{j'}_N\bigl(\bvs^{\ell'}\bigr) = N \xi_{j,j'}\bigl(R^{j,j'}_{\ell, \ell'}\bigr),
\end{equation}
and the Hamiltonian \eqref{eq:hamiltonian} is a centered Gaussian process with covariance
\begin{equation}\label{SHOWLABEL}
\E H_N(\bvs^\ell) H_N(\bvs^{\ell'}) = N \Sum(\bxi(\bR_{\ell,\ell'})),
\end{equation}
where the sum of all entries in a matrix is denoted by 
\begin{equation}\label{SHOWLABEL}
\Sum(\bA) =  \sum_{j,j' \leq n} A_{j,j'}.
\end{equation}

\subsection{The Limit of the Free Energy}

We now define the constrained free energy. Let $\bm{Q} = \big(Q^{j,j'}\big)_{j,j' \leq n}$ be a $n \times n$ symmetric positive semidefinite matrix with off-diagonals, $Q^{j,j'} \in [-1,1]$ for $j \neq j'$ and diagonals $Q^{j,j} = 1$. Given $\epsilon > 0$, we denote the set of spins with constrained self overlaps by
\begin{align}\label{SHOWLABEL}
Q^\epsilon_N &=  \big\{ \bvs \in S_N^n \mmm \| \bR(\bvs ,\bvs) - \bQ \|_\infty \leq \epsilon \big\},
\end{align} 
where $\| \cdot \|_\infty$ is the infinity norm on $n \times n$ matrices. For an external field $\veh = \big(\vh(j)\big)_{j \leq n} \in \R^n$, we define the free energy as
\begin{equation}\label{eq:freeenergy}
F^\epsilon_N(\bvb,\bQ) = \frac{1}{N} \E \log \int_{Q_N^\epsilon} \exp \Big( H_N(\bvs) + \sum_{j \leq n} \vh(j) \sum_{i \leq N} \vs_i(j) \Big) \, d \lambda_N^n(\bvs),
\end{equation}
where the reference measure $\lambda_N^n = \lambda_N^{\otimes n}$ is the product of normalized uniform measures $\lambda_N$ on $S_N$.

We will prove the limit of $\eqref{eq:freeenergy}$ can be expressed as a Parisi type functional. We begin by introducing some notation. Let
\begin{equation}\label{SHOWLABEL}
\Gamma_n = \Big\{\bA \mmm \bA \text{ is a $n \times n$ positive-semidefinite matrix} \Big\},
\end{equation}
denote the space of $n\times n$ matrices, and let
\begin{equation}\label{SHOWLABEL}
\Pi = \Big\{ \pi: [0,1] \to \Gamma_n \mmm \pi \text{ is left-continuous, $\pi(x_1) \leq \pi(x_2)$ for $x_1 \leq x_2$} \Big\}
\end{equation}
denote the space of left-continuous monotone paths on $\Gamma_n$.
The notation $\pi(x_1) \leq \pi(x_2)$ means $\pi(x_2) - \pi(x_1) \in \Gamma_n$. Distances between paths are given by the metric
\begin{equation}\label{eq:metricpathdist}
\mathrm{d}(\pi, \tilde{\pi}) = \int_0^1 \| \pi(x) - \tilde \pi (x) \|_1 \, dx,
\end{equation}
where $\| \bA \|_1 = \sum_{j,j'} |A_{j,j'}|$. These paths are the functional order parameters of $p$-spin models with vector spins.

Consider a discrete path $\pi \in \Pi$ connecting $\bm{0}$ and $\bQ$,
\begin{equation}\label{SHOWLABEL}
\pi(x) = \bQ_k \text{ for } x_{k - 1} < x \leq x_{k} \text{ for $0 \leq k \leq r$},\quad \pi(0) = \bm{0}, \quad \pi(1) = \bQ.
\end{equation}
This path can be encoded with a sequence of real numbers
\begin{align}\label{eq:xseq}
0 = x_{-1} < x_0 < \dots < x_r = 1,
\end{align}
and a monotone sequence of $n \times n$ symmetric positive semi-definite matrices
\begin{align}\label{eq:qseq}
\bzero = \bQ_0 \leq \bQ_1 \leq \dots \leq \bQ_{r} = \bQ.
\end{align}
Recall definition \eqref{eq:ximatrix}, and denote
\begin{equation}\label{eq:theta}
\btheta(\bA) = \big(\theta_{j,j'}(A_{j,j'})\big)_{j,j' \leq n} = \bA \odot \bxi'(\bA) - \bxi(\bA), 
\end{equation}
where $\bxi'(\bA) = (\xi'_{j,j'} (A_{j,j'}))_{j,j' \leq n}$ is the matrix of entry wise derivatives of $\bxi$. The matrix given by
\begin{equation}\label{eq:delta}
\bD_k = \bxi'( \bQ_{k} ) - \bxi'( \bQ_{k-1} ), \quad 1 \leq k \leq r,
\end{equation}
is positive semidefinite. This can seen by applying the Schur product theorem to the Hadamard product representation \eqref{eq:ximatrix}. 

Given a symmetric positive definite matrix $\bL$, for $k \leq r$ we define recursively
\begin{equation}\label{eq:lambda_def}
\bL_{r} = \bL, \quad \bL_{k} = \bL_{k+1} - x_{k} \bD_{k+1} \text{ for $0 \leq k \leq r-1$}.
\end{equation}
Let $| \cdot |$ be the determinant of $n \times n$ matrices and consider the set 
\begin{equation}\label{eq:lagrange_mult_set}
\cL := \cL(\pi) = \{ \bL \in \Gamma_n \mmm |\bL_0| > 0 \}.
\end{equation}
For $\bL \in \cL$ and discrete $\pi \in \Pi$, we define the following functional
\begin{align}
\sP_{\bvb,\bQ}(\bL,\pi) &= \frac{1}{2}\Big[ \tr(\bL \bQ) - n - \log |\bL| + (\bL_0^{-1} \veh, \veh)  + \sum_{0 \leq k \leq r-1} \frac{1}{x_k} \log \frac{|\bL_{k + 1}|}{|\bL_k|} \label{eq:parisifirst} \\
&\quad  - \sum_{0 \leq k \leq r-1} x_k \cdot \Sum \big( \btheta (\bQ_{k + 1}) - \btheta (\bQ_{k})\big)\Big]. \label{SHOWLABEL}
\end{align}
The dependence on $\bvb$ is through the functions $\bxi$ and $\btheta$ defined in $\eqref{eq:ximatrix}$ and $\eqref{eq:theta}$. The following is the main result:
\begin{theo}\label{theo:MAIN}
	For $n \geq 1$ and $\veh \in \R^n$, the limit of the free energy at inverse temperature $\bvb$ and constraint $\bQ$ is given by
	\begin{equation}\label{SHOWLABEL}
	\lim_{\epsilon \to 0}\lim_{N \to \infty} F^\epsilon_N(\bvb, \bQ) = \inf_{\pi,\bL} \sP_{\bvb,\bQ}(\bL,\pi).
	\end{equation}
	The infimum is over $\bL \in \cL$ and discrete paths given by \eqref{eq:xseq} and \eqref{eq:qseq} over all $r \geq 1$.
\end{theo}

\noindent \textbf{Remark.} If $\det(\bQ) = 0$, we show in Lemma~\ref{lem:minimizer} that for all fixed $\beta$ and $\pi$, 
\begin{equation*}
\inf_{\bL} \sP_{\bvb,\bQ}(\bL,\pi) = -\infty.
\end{equation*}
By concentration of measure, this implies that degenerate configurations have exponentially low probability of appearing in the product Gibbs measure.
\\\\
\noindent \textbf{Remark.} Our form of the Parisi functional $\sP_{\bvb,\bQ}(\bL,\pi)$, is missing the $\frac{1}{2}\tr(\bL_0^{-1} \bm\Delta_1)$ that appears in  \cite{PTSPHERE}. This is because we assumed $x_0 > 0$ in \eqref{eq:xseq} while $x_0 = 0$ in  \cite{PTSPHERE}. By applying L'H\^opital's rule and Jacobi's formula, this term can be recovered by observing
\begin{align*}
\lim_{x_0 \to 0} \frac{1}{2x_0} \log \frac{|\bL_1|}{|\bL_1 - x_0 \bD_1|} &=\lim_{x_0 \to 0} \frac{1}{2} |\bL_1 - x_0 \bD_1|^{-1} \tr( |\bL_1 - x_0 \bD_1| (\bL_1 - x_0 \bD_1)^{-1} \bD_1) 
\\&= \frac{1}{2} \tr(\bL_0^{-1} \bD_1).
\end{align*}

\subsection{Outline of the Paper:} We begin by using an analogue of Guerra's interpolation to prove the upper bound in Section~\ref{sec:upperbound}. In Section~\ref{sec:func_RPC}, we prove the sharpness of functionals that appeared in the upper bound using classical large deviations. We begin the proof of the lower bound by using the Poincar\'e limit to derive an analogue of the Aizenman--Sims--Starr scheme for high dimensional spherical spin glass models in Section~\ref{sec:ASS}. In Section~\ref{sec:modpert}, we introduce a perturbation of the Hamiltonian that will force the overlaps under the asymptotic Gibbs measure to satisfy the synchronization properties used in the study of vector spin glass models. In Section~\ref{sec:lwbd} we combine all the results and finish the proof of the lower bound using standard cavity computations.

\section{Upper Bound --- Guerra's Interpolation}\label{sec:upperbound}

We begin by proving the upper bound of the free energy.
\begin{lem} For $n \geq 1$ and $\veh \in \R^n$, 
	\begin{equation}\label{SHOWLABEL}
	\lim_{\epsilon \to 0}\limsup_{N \to \infty} F^\epsilon_N(\bvb, \bQ) \leq \inf_{\bL,\pi} \sP_{\bvb,\bQ}(\bL,\pi).
	\end{equation}
\end{lem}
\noindent A version of this upper bound was proved in Section 2 of \cite{PTSPHERE}. We will provide a different proof using the Ruelle probability cascades and Guerra's interpolation. The main difference is the following proof will hold without the condition that the diagonals of $\bL$ are greater than $1$.

Consider the sequence of real numbers
\begin{equation}\label{eq:rpc1}
0 = x_{-1} < x_0 < \dots < x_r = 1,
\end{equation}
and the sequence of $n \times n$ positive semi definite matrices
\begin{equation}\label{eq:rpc2}
\bm{0} = \bQ_{0}  \leq \bQ_{1} \leq \dots \leq \bQ_r = \bQ.
\end{equation}
Let $(v_\alpha)_{\alpha \in \N^r}$ be the weights of the Ruelle probability cascades \cite{Ruelle} corresponding to the sequence \eqref{eq:rpc1}. For paths $\alpha^1, \alpha^2 \in \N^r$, we denote the common vertices by
\begin{equation}\label{SHOWLABEL}
\alpha^1 \wedge \alpha^2 = \min\big\{ 0 \leq j \leq r \mmm \alpha_1^1 = \alpha_1^2, \dots, \alpha_j^1 = \alpha_j^2, \alpha^1_{j + 1} \neq \alpha^2_{j + 1} \big\}
\end{equation} 
and $\alpha^1 \wedge \alpha^2 = r$ if $\alpha^1 = \alpha^2$. Consider independent centered Gaussian processes $Z(\alpha) = \big (Z^j(\alpha) \big)_{j \leq n}$ and $Y(\alpha)$ indexed with $\alpha \in \N^r$ and covariances
\begin{align}
\Cov( Z(\alpha^1),Z(\alpha^2) ) &=  \bxi'(\bQ_{\alpha^1 \wedge \alpha^2}) , \label{eq:cov_Z} \\
\Cov( Y(\alpha^1),Y(\alpha^2) ) &= \Sum \big( \btheta (\bQ_{\alpha^1 \wedge \alpha^2} ) \big).\label{eq:cov_Y}
\end{align}
Let $Z_i(\alpha)$ be an independent copy of $Z(\alpha)$ also independent of $Y(\alpha)$. A Gaussian interpolation argument will bound the free energy with functions of these Gaussian processes. 

\begin{lem} For all $N > 0$, there exists a constant $L$ such that
	\begin{align}
	F^\epsilon_N(\bvb, \bQ) &\leq  \frac{1}{N} \E \log \sum_{\alpha \in \N^r} v_\alpha \int_{Q^\epsilon_N} \exp \Bigl( \sum_{i \leq N} \sum_{j \leq n} \vs_i(j) \big(Z_i^j(\alpha) + \vh(j)  \big) \Bigr) \, d \lambda_N^n(\bvs) \notag
	\\&\quad - \frac{1}{N} \E \log \sum_{\alpha \in \N^r} v_\alpha \exp \sqrt{N}Y(\alpha) + L \epsilon \label{eq:upbd_first_bd}.
	\end{align}
	
\end{lem}

\begin{proof}
	The result follows from Gaussian interpolation. For $0 \leq t \leq 1$, we define the interpolating Hamiltonian
	\begin{equation*}
	H_{t} (\bvs , \alpha) = \sqrt{t} H_N(\bvs) +  \sum_{i \leq N} \sum_{j \leq n} \vs_i(j) \big(\sqrt{1  - t} Z_i^j(\alpha) + \vh(j) \big)+ \sqrt{t} \sqrt{N} Y(\alpha),
	\end{equation*}
	on $S_N^n \times \N^r$. For a given a constraint $\bQ$, we define the interpolating free energy function
	\begin{equation*}
	\varphi(t) = \frac{1}{N} \E \log \sum_{\alpha \in \N^r} v_\alpha \int_{Q_N^\epsilon} \exp H_{t} (\bvs, \alpha) \, d \lambda_N^n(\bvs).
	\end{equation*}
	
	Let $\langle \cdot \rangle_t$ be the average on $Q_N^\epsilon \times \N^r$ with respect to the Gibbs measure
	\begin{equation*}
	G(d \bvs, \alpha) \propto v_\alpha \exp H_{t} (\bvs, \alpha) \, d \lambda_N^n(\bvs).
	\end{equation*}
	A straightforward computation shows 
	\begin{equation*}
	\phi'(t) = \frac{1}{N} \E \Big\langle \frac{\partial }{\partial t} H_{t} (\bvs, \alpha) \Big\rangle_{t}.
	\end{equation*}
	By Gaussian integration by parts \cite[Lemma 1.1]{PBook},
	\begin{align}
	\frac{1}{N} \E \Big\langle \frac{\partial }{\partial t} H_{t} (\bvs, \alpha) \Big\rangle_{t} &= \frac{1}{2} \E \Big \langle \Sum  \Big( \bxi(\bR_{1,1}) - \bR_{1,1} \odot \bxi'(\bQ_{\alpha^1 \wedge \alpha^1} )  + \btheta(\bQ_{\alpha^1 \wedge \alpha^1}) \Big) \Big \rangle_t  \label{eq:guerra_first_term}
	\\&\quad - \frac{1}{2} \E \Big \langle \Sum  \Big( \bxi(\bR_{1,2}) - \bR_{1,2} \odot \bxi'(\bQ_{\alpha^1 \wedge \alpha^2} )  + \btheta(\bQ_{\alpha^1 \wedge \alpha^2}) \Big) \Big \rangle_t \label{eq:guerra_second_term}.
	\end{align} 
	
	We use convexity to bound \eqref{eq:guerra_second_term}. Since $\vb_p = \vec{0}$ for odd $p$, $\xi_{j,j'}(x)$ is a convex function for all $j,j' \leq n$ and therefore lies above all its tangent lines. That is,
	\begin{equation*}
	\xi_{j,j'}(a) - a \xi'_{j,j'}(b) + \theta(b) \geq 0 \text{ for all $a,b \in \R$}.
	\end{equation*}
	which implies, $\Sum  \big( \bxi(\bR_{1,2}) - \bR_{1,2} \odot \bxi'(\bQ_{\alpha^1 \wedge \alpha^2} )  + \btheta(\bQ_{\alpha^1 \wedge \alpha^2}) \big)$ is non-negative. To bound \eqref{eq:guerra_first_term} we use definition \eqref{eq:theta} and notice \eqref{eq:guerra_first_term} is equal to
	\begin{align}
	\E \Big \langle \Sum\Big( \bxi(\bR_{1,1}) - \bxi(\bQ_{\alpha^1 \wedge \alpha^1}) - (\bR_{1,1} - \bQ_{\alpha^1 \wedge \alpha^1})\odot \bxi'(\bQ_{\alpha^1 \wedge \alpha^1}) \Big)  \Big \rangle_t \label{eq:upbd2}.
	\end{align} 
	The self overlaps are constrained, so $\| \bR_{1,1} - \bQ_{\alpha^1 \wedge \alpha^1} \|_\infty \leq \epsilon$. Continuity of $\bxi$ implies \eqref{eq:guerra_first_term} is bounded by $L \epsilon$, for some constant $L$ that does not depend on $N$.
	
	These bounds on \eqref{eq:guerra_first_term} and \eqref{eq:guerra_second_term} imply
	\begin{equation}\label{eq:guerraupbd}
	\phi'(t) \leq L \epsilon.
	\end{equation}
	By the mean value theorem, \eqref{eq:guerraupbd} gives us the upper bound
	\begin{equation}\label{eq:upbd_ineq}
	\phi(1) \leq \phi(0) + L\epsilon,
	\end{equation}
	where
	\begin{align}
	\phi(1) &= F^\epsilon_N(\bvb, \bQ)  + \frac{1}{N} \E \log \sum_{\alpha \in \N^r} v_\alpha \exp \sqrt{N}Y(\alpha), \\
	\phi(0) &= \frac{1}{N} \E \log \sum_{\alpha \in \N^r} v_\alpha \int_{Q^\epsilon_N} \exp \Big( \sum_{i \leq N} \sum_{j \leq n} \vs_i(j) \big(Z_i^j(\alpha) + \vh(j)  \big) \Big) \, d \lambda_N^n(\bvs) \label{eq:z_func}.
	\end{align}
	Rearranging terms finishes the proof of the upper bound.
\end{proof}

The terms in \eqref{eq:upbd_first_bd} containing $Y(\alpha)$ and $Z(\alpha)$ can be computed explicitly using the recursive construction of the Ruelle probability cascades \cite[Theorem 2.9]{PBook}. Recalling the covariance structure in \eqref{eq:cov_Y}, a recursive computation \cite[Chapter 3]{PBook} shows
\begin{align}
\limsup_{N \to \infty} \frac{1}{N} \E \log \sum_{\alpha \in \N^r} v_\alpha \exp \sqrt{N}Y(\alpha) = \sum_{0 \leq k \leq r-1} x_k \cdot \Sum \big( \btheta (\bQ_{k + 1}) - \btheta (\bQ_{k})\big) \label{eq:upbd_Y_recur}.
\end{align}

The term in \eqref{eq:upbd_first_bd} containing $Z(\alpha)$ can be computed similarly after decoupling the constraint on $Q_N^\epsilon$ using Lagrange multipliers and rotational invariance \cite[Lemma 1]{PTSPHERE}. Let $\nu_{N}$ be the standard Gaussian measure on $\R^N$. We write $\bo(j)\in \R^N$ in its polar coordinate form $\bo(j) = (s_j \bs(j))$, where $s_j = \frac{\|\bo(j)\|}{\sqrt{N}} \in \R^+$ and $\bs(j) =  \frac{\sqrt{N} \bo(j)}{\|\bo(j)\|} \in S_N$. Let $\gamma_N$ denote the law of $s_j$ under $\nu_N$. By rotational invariance, the law of $\bs(j)$ under $\nu_j$ is $\lambda_N$, and $\bs(j)$ and $s_j$ are independent. We express \eqref{eq:z_func} in terms of a Gaussian integral.

\begin{lem} \label{lem:gaussian_int}
	There exists a $\delta \in (0,\epsilon)$, such that \eqref{eq:z_func} is bounded above by
	\begin{align}
	\frac{1}{N} \E \log \sum_{\alpha \in \N^r} v_\alpha \int_{\Omega^{\epsilon, \delta}_N} \exp \Big( \sum_{i \leq N} \sum_{j \leq n} \vo_i(j) \big(Z_i^j(\alpha) + \vh(j)  \big) \Big) \, d \nu_N^n(\bvo)  - \frac{n \log \nu_N( E^\delta_N ) }{ N } + L \delta \label{eq:surface_to_gaussian}
	\end{align}
	where the $\delta$ shell around $Q_N^\epsilon$ is denoted by
	\begin{equation}\label{eq:deltashell2}
	\Omega_N^{\epsilon,\delta} = \big\{  \bvo = (s_j \bs(j))_{j \leq n} \in (\R^N)^n \mmm \bvs \in Q_N^{\epsilon}, \, s_j \in [\sqrt{1 - \delta}, \sqrt{1 + \delta}] \text{ for all } j \leq n \big\}
	\end{equation}
	and the $\delta$ neighbourhood of the radial component is denoted by
	\begin{equation*}
	E_N^\delta = \{ x \in \R^N \mmm \|x\| \in [ \sqrt{(1 - \delta) N}, \sqrt{(1 + \delta)N} ] \}.
	\end{equation*}
\end{lem}

\begin{proof}
	We will use a Gaussian interpolation argument. Let $\tilde{Z}^j_i(\alpha)$ be an independent copy of $Z_i^j(\alpha)$. For $0 \leq t \leq 1$, we define the interpolating Hamiltonian
	\begin{equation*}
	H_{t} (\bvo , \alpha) = \sqrt{t} \Big( \sum_{i \leq N} \sum_{j \leq n} \vs_i(j) \tilde{Z}_i^j(\alpha) \Big) + \sqrt{1 - t} \Big( \sum_{i \leq N} \sum_{j \leq n} \vo_i(j)Z_i^j(\alpha)\Big) + \sum_{i \leq N} \sum_{j \leq n} \vs_i(j) \vh(j) ,
	\end{equation*}
	on $\Omega_N^{\epsilon,\delta} \times \N^r$. The corresponding interpolating free energy function is denoted by
	\begin{equation*}
	\varphi(t) = \frac{1}{N} \E \log \sum_{\alpha \in \N^r} v_\alpha \int_{\Omega_N^{\epsilon,\delta}} \exp H_{t} (\bvo, \alpha) \, d \nu_N^n(\bvo).
	\end{equation*}
	Let $\langle \cdot \rangle_t$ be the average on $\Omega_N^{\epsilon,\delta} \times \N^r$ with respect to the Gibbs measure
	\begin{equation*}
	G(d \bvo, \alpha) \propto v_\alpha \exp H_{t} (\bvo, \alpha) \, d \nu_N^n(\bvo).
	\end{equation*}
	By Gaussian integration by parts,
	\begin{equation*}
	\phi'(t) = \frac{1}{N} \E \Big\langle \frac{\partial }{\partial t} H_{t} (\bvo, \alpha) \Big\rangle_{t} = \E \bigg\langle \E \frac{\partial H_t(\bvo^1, \alpha^1) }{ \partial t} \cdot H_t(\bvo^1, \alpha^1) - \E \frac{\partial 	H_t(\bvo^1, \alpha^1)}{ \partial t} \cdot 	H_t(\bvo^2, \alpha^2) \bigg\rangle_t.
	\end{equation*}
	Computing the covariances, we get
	\begin{align*}
	\phi'(t) &= \frac{1}{2} \E \Big \langle \Sum  \Big(\bR(\bvs^1,\bvs^1) \odot \bxi'(\bQ_{\alpha^1 \wedge \alpha^1}) - \bR(\bvo^1,\bvo^1) \odot \bxi'(\bQ_{\alpha^1 \wedge \alpha^1})  \Big) \Big \rangle_t  
	\\&\quad - \frac{1}{2} \E \Big \langle \Sum  \Big( \bR(\bvs^1,\bvs^2)\odot \bxi'(\bQ_{\alpha^1 \wedge \alpha^2}) - \bR(\bvo^1,\bvo^2) \odot \bxi'(\bQ_{\alpha^1 \wedge \alpha^2})  \Big) \Big \rangle_t.
	\end{align*} 
	Since $\vo_i(j) = s_j \vs_i(j)$ and $s_j \in [\sqrt{1 - \delta}, \sqrt{1 + \delta}]$, we have the bound
	\begin{equation*}
	\Sum  \Big( \bR(\bvs^1,\bvs^2)\odot \bxi'(\bQ_{\alpha^1 \wedge \alpha^2}) - \bR(\bvo^1,\bvo^2) \odot \bxi'(\bQ_{\alpha^1 \wedge \alpha^2})  \Big) \leq \delta  n^2  \| \bxi'(1) \|_\infty.
	\end{equation*}
	By the triangle inequality,
	\begin{equation*}
	|\phi'(t) | \leq n^2 \delta \| \bxi'(1) \|_\infty= L \delta,
	\end{equation*}
	resulting in the bound
	\begin{equation}\label{eq:upbd_second_interp}
	\phi(1) \leq \phi(0) + L \delta.
	\end{equation}
	The ending term of the interpolation can be simplified using rotational invariance of $\nu_N$, 
	\begin{align}
	\phi(1) &= \frac{1}{N} \E \log \sum_{\alpha \in \N^r} v_\alpha \int_{\Omega_N^{\epsilon,\delta}} \exp  \biggl( \sum_{i \leq N} \sum_{j \leq n} \vs_i(j) \big(\tilde{Z}_i^j(\alpha) + \vh(j) \big) \biggr) \, d \nu_N^n(\bvo) \notag
	\\&= \frac{1}{N} \E \log \sum_{\alpha \in \N^r} v_\alpha \int_{ [\sqrt{1 - \delta}, \sqrt{1 + \delta} ]^n } \int_{Q_N^{\epsilon}} \exp  \biggl( \sum_{i \leq N} \sum_{j \leq n} \vs_i(j) \big(\tilde{Z}_i^j(\alpha) + \vh(j) \big) \biggr) \, d \lambda_N^n(\bvs) d \gamma^n_N(s) \notag
	\\&= \frac{1}{N} \E \log \sum_{\alpha \in \N^r} v_\alpha \int_{Q_N^{\epsilon}} \exp  \biggl( \sum_{i \leq N} \sum_{j \leq n} \vs_i(j) \big(\tilde{Z}_i^j(\alpha) + \vh(j) \big) \biggr) \, d \lambda_N^n(\bvs) + \frac{n \log \nu_N( E_N^\delta ) }{ N }\label{eq:upbd_second_interp2}.
	\end{align}
	Substituting \eqref{eq:upbd_second_interp2} into \eqref{eq:upbd_second_interp} gives the bound
	\begin{align}
	&\frac{1}{N} \E \log \sum_{\alpha \in \N^r} v_\alpha \int_{\Omega_N^{\epsilon,\delta}} \exp  \biggl( \sum_{i \leq N} \sum_{j \leq n} \vs_i(j) \big(\tilde{Z}_i^j(\alpha) + \vh(j) \big) \biggr) \, d \nu_N^n(\bvo) \label{eq:upbd_third_step}
	\\&\leq \frac{1}{N} \E \log \sum_{\alpha \in \N^r} v_\alpha \int_{\Omega_N^{\epsilon,\delta}} \exp  \biggl( \sum_{i \leq N} \sum_{j \leq n} \Big( s_j \vs_i(j) \tilde{Z}_i^j(\alpha) + \vs_i(j) \vh(j) \Big)  \biggr) \, d \nu_N^n(\bvo) - \frac{n \log \nu_N( E_N^\delta ) }{ N } + L \delta.\notag
	\end{align}
	On the set $\Omega_N^{\epsilon,\delta}$, the Cauchy--Schwarz inequality implies 
	\begin{equation*}
	\bigg| \sum_{j \leq n} \sum_{i \leq N} \Big(s_j \vs_i(j) \vh(j) - \vs_i(j) \vh(j) \Big) \bigg| \leq \delta \sqrt{N} \sum_{j \leq n}\| \bvs(j) \| \cdot |\vh(j)| \leq \delta L N.
	\end{equation*}
	Therefore, we can replace $\vs_i(j) \vh(j)$ with $\vo_i(j) \vh(j)$ in the upper bound of \eqref{eq:upbd_third_step} and absorb the error into $L \delta$ giving
	\begin{equation*}
	\frac{1}{N} \E \log \sum_{\alpha \in \N^r} v_\alpha \int_{\Omega^{\epsilon, \delta}_N} \exp \Big( \sum_{i \leq N} \sum_{j \leq n} \vo_i(j) \big(Z_i^j(\alpha) + \vh(j)  \big) \Big) \, d \nu_N^n(\bvo) - \frac{n \log \nu_N( E^\delta_N ) }{ N } + L \delta,
	\end{equation*}
	the required upper bound in \eqref{eq:surface_to_gaussian}.
\end{proof}

We now explicitly compute the upper bound of \eqref{eq:surface_to_gaussian}. We denote the subset of $\R^{Nn}$ constrained by coupling the overlaps with,
\begin{equation}\label{eq:modeified_Omega}
\tilde{\Omega}_N^\epsilon = \big\{\, \bvo \in (\R^N)^n \mmm R^{j,j'}\bigl(\bvo, \bvo\bigr) \in [Q^{j,j'} - \epsilon, Q^{j,j'} + \epsilon] \text{ for all $j,j' \leq n$} \,\big\}.
\end{equation}
For $\delta < \epsilon$, $\Omega_N^{\epsilon,\delta} \subset \tilde{\Omega}_N^{2\epsilon}$ so \eqref{eq:surface_to_gaussian} is bounded above by
\begin{equation}\label{eq:upbd3}
\frac{1}{N} \E \log \sum_{\alpha \in \N^r} v_\alpha \int_{\tilde\Omega^{2\epsilon}_N} \exp \Big( \sum_{i \leq N} \sum_{j \leq n} \vo_i(j) \big(Z_i^j(\alpha) + \vh(j)  \big) \Big) \, d \nu_N^n(\bvo)  - \frac{n \log \nu_N( E^\delta_N ) }{ N } + L \epsilon.
\end{equation}
For any $\bvo \in \tilde\Omega^{2\epsilon}_N$ and $\bL \in \cL$,
\begin{equation*}
\bigg\| \sum_{j,j' \leq n} \Lambda^{j,j'} Q^{j,j'} - \frac{1}{N} \sum_{j,j' \leq n} \sum_{i \leq N} \Lambda^{j,j'} \vo_i(j') \vo_i(j) \bigg\|_1 \leq 2 \epsilon \| \bL \|_1.
\end{equation*}
Therefore, adding and subtracting $\frac{1}{2} \sum_{i \leq N} ((\bL - \bI) \vo_i, \vo_i)$ from the exponent implies \eqref{eq:upbd3} can be bounded above by
\begin{align}
&\frac{1}{N} \E \log \sum_{\alpha \in \N^r} v_\alpha \int_{\tilde{\Omega}_N^{2 \epsilon}} \exp \Bigl( \sum_{i \leq N} \sum_{j \leq n} \vo_i(j) \big(Z_i^j(\alpha) + \vh(j)  \big) - \frac{1}{2}  \sum_{i \leq N} ((\bL - \bI) \vo_i, \vo_i)  \Bigr) \, d \nu_N^n(\bvo) \notag
\\&+ \frac{1}{2} \tr(\bL \bQ)- \frac{n}{2}  - \frac{n \log \nu_N( E^\delta_N ) }{ N } + 2\epsilon\|\bL\|_1 - L \epsilon \label{eq:upbd4}.
\end{align}
Since $\tilde{\Omega}_N^{2 \epsilon} \subset (\R^N)^n$, if we define the function,
\begin{equation}\label{eq:upbdrec_y}
Y_{r,i}(\alpha) =  \frac{1}{(2\pi)^{n/2}} \int_{\R^n} \exp \Big(  \sum_{j \leq n} \vo_i(j) \big(Z_i^j(\alpha) + \vh(j)  \big)  - \frac{1}{2} \sum_{j, j' \leq n} \Lambda^{j,j'} \vo_i(j) \vo_i(j') \Big) \, d \vo_i
\end{equation}
then our upper bound \eqref{eq:upbd4} can be written as
\begin{equation}
\frac{1}{N} \E \log \sum_{\alpha \in \N^r} v_\alpha \prod_{i \leq N} Y_{r,i}(\alpha)  + \frac{1}{2} \tr(\bL \bQ)- \frac{n}{2}  - \frac{n \log \nu_N( E^\delta_N ) }{ N } + \epsilon\|\bL\|_1 - L \epsilon.  \label{eq:rec_upbd_step2} 
\end{equation}

The term containing $Y_{r,i}(\alpha)$ in \eqref{eq:rec_upbd_step2} can be computed recursively. Let $\bvz_k = (z^j_{k})_{j \leq n}$ be a Gaussian vector with covariance $\bD_k$ defined in \eqref{eq:delta} and let $\bvz_k$ be independent for $1 \leq k \leq r$. For $i \leq M$ let $\bvz_{k,i}$  be an independent copy of $\bvz_k$. We define the recursion starting with
\begin{equation}\label{eq:recur}
Y_{r,i} = \log \frac{1}{(2\pi)^{n/2}} \int_{\R^n} \exp \bigg(  \sum_{j \leq n} \vo_i(j) \Big(\sum_{1 \leq k \leq r} z_{k,i}^j + \vh(j)  \Big)  - \frac{1}{2} \sum_{j, j' \leq n} \Lambda^{j,j'} \vo_i(j) \vo_i(j') \bigg) \, d \vo_i
\end{equation}
with subsequent values for $0 \leq k \leq r - 1$ given recursively by
\begin{equation}\label{eq:recursion}
Y_{k,i} = \frac{1}{x_k} \log \E_k \exp x_k Y_{k + 1,i},
\end{equation}
where $\E_k$ refers to expectation with respect to the random vector $\bvz_{k + 1,i}$. The $\bvz_i$ are \iid so $Y_{0,i} = Y_{0,1}$ for all $i \leq N$. The recursive representation of the average in \cite[Theorem 2.9]{PBook} implies \eqref{eq:rec_upbd_step2}  can be written as
\begin{equation}\label{eq:rec_upbd_step3}
Y_{0,1}   + \frac{1}{2} \tr(\bL \bQ)- \frac{n}{2}  - \frac{n \log \nu_N( E^\delta_N ) }{ N } + \epsilon\|\bL\|_1 - L \epsilon.
\end{equation} 

In this model, $Y_{0,1}$ has a closed form. Starting from the start of the recursion, a direct computation (see equation (2.17) in \cite{PTSPHERE}) shows
\begin{equation}\label{eq:lim_rec_4}
Y_{r,1} = -\frac{1}{2} \log |\bL| + \frac{1}{2} \bigg(  \bL^{-1} \Big( \sum_{1 \leq k \leq r} \vec{z}_{k,1} + \veh \Big), \Big( \sum_{1 \leq k \leq r} \vec{z}_{k,1} + \veh \Big)  \bigg).
\end{equation}
Here $(\cdot, \cdot )$ is the scalar product of vectors in $\R^n$. The first term is non-random and will propagate through the recursion. The second term can be computed recursively using the following result:
\begin{lem}\label{lem:calc_lem_1}
	Let $\vg$ be a Gaussian vector with covariance $\bC$. Then for any $y \in \R^n$ and $x \in (0,1]$, 
	\begin{align*}
	\frac{1}{x} \log \E \exp \Big( \frac{x}{2} \big( \bA^{-1} (y + \vg\,), y + \vg \,\big) \Big) = \frac{1}{2x} \log \frac{|\bA|}{|\bA - x\bC|} + \frac{1}{2} \big( (\bA - x\bC)^{-1} y , y \big).
	\end{align*}
\end{lem}
\begin{proof}
	The one dimensional case was proven in \cite[Lemma 3.5]{TSPHERE}. We will prove the analogous result for $\R^n$. The expectation can be computed explicitly as follows,
	\begin{align*}
	& \E \exp \Big( \frac{x}{2} \big( \bA^{-1} (y + \vg\,), y + \vg\, \big) \Big) 
	\\&= \bigg(\frac{|\bC|^{-1}}{ (2 \pi)^{n}}\bigg)^{1/2}\int_{ \R^n } \exp \Big( \frac{x}{2} \big( \bA^{-1} (y + z\,), y + z \, \big) - \frac{1}{2} \big( \bC^{-1} z, z \,\big) \Big) \, d z 
	\\&= \bigg(\frac{|\bC|^{-1}}{(2 \pi)^{n}}\bigg)^{1/2}\int_{ \R^n } \exp \Big( \frac{x}{2} \big( (\bA - x\bC)^{-1} y, y \big) - \frac{1}{2} \big( (\bC^{-1} - x\bA^{-1}) (z - \bB y \,), (z - \bB y \,) \big) \Big) \, d z 
	\\&= \bigg( \frac{|\bC|^{-1}}{|\bC^{-1} - x\bA^{-1}|} \bigg)^{1/2} \exp \Big( \frac{x}{2} \big( (\bA - x\bC)^{-1} y, y \,\big) \Big)
	\end{align*}
	where the matrix $\bB$ is given by
	\begin{equation*}
	\bB =  x(\bC^{-1} - x \bA^{-1})^{-1}\bA^{-1}.
	\end{equation*}
	The conclusion follows immediately if we rewrite the matrices in the normalizing constant as,
	\begin{equation*}
	(\bC^{-1} - x \bA^{-1}) = \bC^{-1}(\bA - x\bC) \bA^{-1}, 
	\end{equation*}
	which implies
	\begin{equation*}
	|\bC^{-1} - x \bA^{-1}| = |\bC|^{-1} |\bA - x \bC| |\bA|^{-1}.
	\end{equation*}
\end{proof}

Using Lemma~\ref{lem:calc_lem_1} to compute the recursion gives the appropriate closed form.
\begin{cor}\label{cor:closed_form_recur}
	If $|\bL_0| > 0$, then
	\begin{align}
	Y_{0,1} = -\frac{1}{2} \log |\bL| + \frac{1}{2} \Big( \bL_0^{-1} \vec{h}, \vec{h} \Big) + \frac{1}{2} \sum_{0 \leq k \leq r-1} \frac{1}{x_k}  \log \frac{|\bL_{k + 1}|}{|\bL_k|}\label{eq:closedformY0}.
	\end{align}
\end{cor}

\begin{proof}
	Using Lemma \ref{lem:calc_lem_1} to compute the expectation of the second term in \eqref{eq:lim_rec_4} recursively implies
	\begin{align}
	Y_{r-1,1} &= -\frac{1}{2} \log |\bL| + \frac{1}{2x_{r}} \log \frac{|\bL_{r}|}{|\bL_{r} - x_{r-1} \bD_r |} \label{eq:rec_y_1}
	\\&\quad + \frac{1}{2} \bigg( (\bL_{r} - x_{r-1} \bD_r)^{-1} \Big( \sum_{1 \leq k \leq r -1} \vec{z}_{k,1} + \vh \,\Big), \Big( \sum_{1\leq k \leq r - 1} \vec{z}_{k,1} + \vh \, \Big) \bigg) \label{eq:rec_y_2}.
	\end{align}
	Again, the terms in \eqref{eq:rec_y_1} are non-random, so they propagate through the recursion. Computing the terms in \eqref{eq:rec_y_2} inductively using repeated applications of Lemma \ref{lem:calc_lem_1} implies
	\begin{align}\label{eq:upbd_2}
	Y_{0,1} &= -\frac{1}{2} \log |\bL| + \frac{1}{2} \Big( \bL_0^{-1} \vec{h}, \vec{h} \Big) + \frac{1}{2} \sum_{0 \leq k \leq r-1} \frac{1}{x_k}  \log \frac{|\bL_{k + 1}|}{|\bL_k|}.
	\end{align}
\end{proof}

Notice $\nu_N$ concentrates around the sphere of radius $\sqrt{N}$ in high dimensions, so the Gaussian term will vanish in the limit by classical large deviations \cite[Lemma 3.2 and Lemma 3.3]{TSPHERE}. That is, for any $\delta > 0$,
\begin{equation}\label{eq:error_large_dev}
\limsup_{N \to \infty} \frac{n \log \nu_N( E^\delta_N ) }{ N } = 0.
\end{equation}
The other terms vanish by taking $\epsilon \to 0$, so combining \eqref{eq:closedformY0} and \eqref{eq:error_large_dev} with \eqref{eq:rec_upbd_step3}, gives the bound
\begin{align}
& \lim_{\epsilon \to 0}\limsup_{N \to \infty} \frac{1}{N} \E \log \sum_{\alpha \in \N^r} v_\alpha \int_{Q^\epsilon_N} \exp \Big( \sum_{i \leq N} \sum_{j \leq n} \vs_i(j) \big(Z_i^j(\alpha) + \vh(j)  \big) \Big) \, d \lambda_N^n(\bvs) \notag
\\&\leq \frac{1}{2} \bigg( \tr(\bL \bQ) - n - \log |\bL|  + (\bL^{-1}_0 \veh,\veh) +  \sum_{0 \leq k \leq r-1} \frac{1}{x_k}  \log \frac{|\bL_{k + 1}|}{|\bL_k|} \bigg).\label{eq:upbd_Z_recur}
\end{align}

The upper bound in \eqref{eq:upbd_Z_recur} holds for all $\bL \in \cL$. Applying the bounds \eqref{eq:upbd_Z_recur} and \eqref{eq:upbd_Y_recur} to \eqref{eq:upbd_ineq} and taking the infimum over all discrete paths encoded by the monotone sequences \eqref{eq:rpc1} and \eqref{eq:rpc2} shows
\begin{equation}\label{eq:upbound}
\lim_{\epsilon \to 0}\limsup_{N \to \infty} F^\epsilon_N(\bvb, \bQ) \leq \inf_{\bL,\pi} \sP_{\bvb,\bQ}(\bL,\pi),
\end{equation}
completing the proof of the upper bound. 

	\section{Sharpness of the Upper Bound}\label{sec:func_RPC}

We now prove for every fixed path $\pi$, the upper bound \eqref{eq:upbd_Z_recur} is asymptotically sharp in the sense that it attains equality after minimizing over $\bL$. This fact will be used again when a similar functional appears in the proof of the lower bound.
The proof of this sharpness for the replica symmetric case can be found in \cite[Lemma 4]{PTSPHERE}. We will provide a proof of the general  case below.

Let $\pi$ be any fixed discrete monotone path characterized by the sequences \eqref{eq:rpc1} and \eqref{eq:rpc2} and denote the functional appearing in \eqref{eq:surface_to_gaussian} by
\begin{equation}\label{SHOWLABEL}
f_N^1(\pi) = \frac{1}{N} \E \log  \sum_{\alpha \in \N^r} v_\alpha \int_{\Omega_N^{\epsilon,\delta}} \exp \Big( \sum_{i \leq N} \sum_{j \leq n} \vo_i(j) \big( Z_i^j(\alpha) + \vh(j)\big) \Big) \, d \nu_{N}^n(\bvs).
\end{equation}
We will prove the matching lower bound of \eqref{eq:upbd_Z_recur} by decoupling the functional $f_N^1(\pi)$ from the constraint $\bQ$ and explicitly computing its value recursively.
\begin{lem}\label{lem:sharpupbd} For all $0 < \delta < \epsilon$,
	\begin{equation*}
	\liminf_{N \to \infty} f_N^1(\pi) \geq \inf_{\bL} \frac{1}{2} \bigg( \tr(\bL \bQ) - n - \log |\bL|  + (\bL^{-1}_0 \veh,\veh) +  \sum_{0 \leq k \leq r-1} \frac{1}{x_k}  \log \frac{|\bL_{k + 1}|}{|\bL_k|} \bigg).
	\end{equation*}
\end{lem}

Recall \eqref{eq:modeified_Omega}, the subset of $\R^{Nn}$ constrained by coupling the overlaps,
\begin{equation}\label{eq:V0}
\tilde{\Omega}_N^\delta = \big\{\, \bvo \in (\R^N)^n \mmm \| \bR\bigl(\bvo, \bvo\bigr) - \bQ \|_\infty \leq \delta \,\big\}.
\end{equation}
Clearly, there exists a $\delta^* < \epsilon$ such that $\Omega_N^{\epsilon,\delta}\supseteq \tilde{\Omega}_N^{\delta^*}$, so
\begin{equation}\label{eq:step1.1}
f_N^1(\pi) \geq \frac{1}{N} \E \log \sum_{\alpha \in \N^r} v_\alpha \int_{\tilde{\Omega}_N^{\delta^*}} \exp \Big( \sum_{i \leq N} \sum_{j \leq n} \vo_i(j) \big(Z_i^j(\alpha) + \vh(j)  \big) \Big) \, d \nu_N^n(\bvo).
\end{equation}

We introduce the Lagrange multipliers $\bL \in \cL$ defined in \eqref{eq:lagrange_mult_set}. Like in the proof of the upper bound, since $\|\bR(\bvo,\bvo) - \bQ\|_\infty \leq \delta^*$ for $\bvo \in \tilde{\Omega}_N^{\delta^*}$, adding and subtracting the quadratic form $\frac{1}{2} \sum_{i \leq N} ((\bL - \bI) \vo_i, \vo_i)$ from the exponent implies \eqref{eq:step1.1} is bounded below by
\begin{align}
&\frac{1}{N} \E \log \sum_{\alpha \in \N^r} v_\alpha \int_{\tilde{\Omega}_N^{\delta^*}} \exp \Bigl( \sum_{i \leq N} \sum_{j \leq n} \vo_i(j) \big(Z_i^j(\alpha) + \vh(j)  \big) - \frac{1}{2}  \sum_{i \leq N} ((\bL - \bI) \vo_i, \vo_i)  \Bigr) \, d \nu_N^n(\bvo)\notag
\\&+ \frac{1}{2} \tr(\bL \bQ)- \frac{n}{2} - \delta^* \|\bL\|_1  \label{eq:boundapprox}.
\end{align}

We view the quantity on the first line of \eqref{eq:boundapprox} as a function of $\bL$ and the region of integration. In general, we denote this integral over sets $V \subset (\R^N)^n$ by
\begin{equation}
\!\!\!\!\!\!\!\!\!\!\!\!
\Phi_V(\bL) \!=\! \frac{1}{N} \E \log \sum_{\alpha \in \N^r} v_\alpha \int_{V} \! \exp \Bigl( \sum_{i \leq N} \sum_{j \leq n} \vo_i(j) \big(Z_i^j(\alpha) + \vh(j)  \big) - \frac{1}{2}  \sum_{i \leq N} ((\bL - \bI) \vo_i, \vo_i)  \Bigr) \, d \nu_N^n(\bvo) \label{eq:integralsubset}
\end{equation}
and the integral over the whole space by
\begin{equation}\label{SHOWLABEL}
F(\bL) := \Phi_{\R^{Nn}}(\bL).
\end{equation}
The function $F(\bL)$ does not depend on $N$, and was computed using the recursion \eqref{eq:recursion} giving the closed form in Corollary~\ref{cor:closed_form_recur},
\begin{equation}\label{eq:minimizefunction}
F(\bL) = \frac{1}{2} \bigg( - \log |\bL| + \Big( \bL_0^{-1} \veh, \veh \, \Big) +  \sum_{0 \leq k \leq r-1} \frac{1}{x_k}  \log \frac{|\bL_{k + 1} |}{|\bL_k|} \bigg).
\end{equation}
We will prove that minimizing over $\bL$ removes the dependence on the constraint $\bQ$ asymptotically. We start by showing there exists a unique $\bL_*$ that minimizes $\frac{1}{2} \tr(\bL \bQ ) + F(\bL)$ if the lower bound \eqref{eq:lowbd_rpc} is finite.

\begin{lem}\label{lem:minimizer}
	Given a positive semi-definite constraint $\bQ$:
	\begin{enumerate}
		\item If $\bQ$ is degenerate, then
		\begin{equation}\label{eq:degeneartediverge}
		\inf_{\bL} \Big( \frac{1}{2} \tr(\bL \bQ ) + F(\bL) \Big) = -\infty.
		\end{equation}
		\item If $\bQ$ is non-degenerate, then there exists a $\bL_* \in \cL$ that minimizes $\frac{1}{2} \tr(\bL \bQ ) + F(\bL)$ and satisfies
		\begin{equation}\label{eq:minimizer}
		\frac{\partial}{\partial t} \Big( \frac{1}{2} \tr\bigl( (\bL_* + t\bB) \bQ \bigr) + F(\bL_* + t\bB) \Big)\Big|_{t = 0} = 0
		\end{equation}
		for any symmetric matrix $\bB$.
	\end{enumerate}
	
\end{lem}

\begin{proof} Consider the eigendecomposition of $\bL \in \cL$,
	\begin{equation*}
	\bL = \bP \bd \bP^\trans.
	\end{equation*}
	Using this change of variables and \eqref{eq:lambda_def}, we see \eqref{eq:minimizefunction} can be rewritten in terms of $\bP$ and $\bd$ as,
	\begin{align}
	\frac{1}{2} \tr(\bL \bQ ) + F(\bL) &= \frac{1}{2} \bigg( \tr\bigl( \bd \bP^\trans \bQ \bP \bigr)  - \log |\bd| + \Big( \Big( \bd - \sum_{0 \leq k < r} x_k \bP^\trans  \bD_{k + 1} \bP  \Big)^{-1} (\bP^\trans\vh), (\bP^\trans\vh) \, \Big) \notag \\ 
	&\quad +  \sum_{0 \leq k \leq r-1} \frac{1}{x_k}  \log \frac{|\bd - \sum_{k + 1 \leq \ell < r}x_\ell \bP^\trans \bD_{\ell + 1} \bP|}{|\bd - \sum_{k \leq \ell < r}x_{\ell} \bP^\trans \bD_{\ell + 1} \bP |} \bigg) \label{eq:tracebound_transformed}.
	\end{align}
	In this form, the infimum is over positive semidefinite diagonal matrices $\bd$ and orthogonal matrices $\bP$ such that $|\bd - \sum_{0 \leq \ell < r}x_{\ell} \bP^\trans \bD_{\ell + 1} \bP | > 0$,
	\begin{equation*}
	\inf_{\bL} \Big( \frac{1}{2} \tr(\bL \bQ ) + F(\bL) \Big) = \inf_{\bd,\bP} \Big( \frac{1}{2} \tr( \bP \bd \bP^\trans \bQ ) + F(\bP \bd \bP^\trans) \Big).
	\end{equation*}
	\textit{Case $\mathrm{(1)}$:} Suppose $\bQ$ is degenerate, i.e. $|\bQ| = 0$. There exists an orthogonal matrix $\bP$, corresponding to the eigendecomposition of $\bQ$, such that $\tilde{\bd} = \bP^\trans \bQ \bP$ and $\tilde{\bd}_{11} = 0$. Given this $\bP$, we choose diagonal matrix $\bd$ with diagonal entries large enough such that all the Gershgorin discs of $\bd - \sum_{k \leq \ell < r}x_{\ell} \bP^\trans \bD_{\ell + 1} \bP$ are contained in the positive real half plane for all $k \leq r -1$. In particular for all $k \leq r-1$, the smallest eigenvalue of $\bd - \sum_{k \leq \ell < r}x_{\ell} \bP^\trans \bD_{\ell + 1} \bP$ is strictly positive and will remain bounded away from zero if we increase the value of the first diagonal element. That is, there exists a $c > 0$ such that
	\begin{equation}\label{eq:smalleig}
	\liminf_{D_{11} \to \infty} \lambda_{min} \Bigl(\bd - \sum_{k \leq \ell < r}x_{\ell} \bP^\trans \bD_{\ell + 1} \bP \Bigr) \geq c > 0 \text{ for all $k \leq r -1$}.
	\end{equation}
	We fix all entries $\bd_{jj}$ for $2 \leq j \leq n$ and show \eqref{eq:tracebound_transformed} diverges to $-\infty$ as we take the first entry $\bd_{11} \to \infty$. For the above choice of $\bd$ and $\bP$, we have
	\begin{align*}
	\tr\bigl( \bd \bP^\trans \bQ \bP \bigr) - \log |\bd| &= \sum_{i = 1}^n \bd_{ii} (\bP^\trans \bQ \bP)_{ii} - \sum_{i = 1}^n \log \bd_{ii}
	\\&= \sum_{i = 2}^n \bd_{ii} (\bP^\trans \bQ \bP)_{ii} - \sum_{i = 2}^n \log \bd_{ii} - \log \bd_{11},
	\end{align*}
	which implies
	\begin{equation}\label{eq:trivcase_1}
	\lim_{D_{11} \to \infty} \Big( \tr\bigl( \bd \bP^\trans \bQ \bP \bigr) - \log |\bd| \Big) = -\infty.
	\end{equation} 
	
	We will now show that the remaining terms of \eqref{eq:tracebound_transformed} are finite. Let $\nu_{min}^{k}$ denote the smallest eigenvalue of $\bd - \sum_{k \leq \ell < r}x_{\ell} \bP^\trans \bD_{\ell + 1} \bP$. By \eqref{eq:smalleig}, we have $\liminf_{D_{11} \to \infty} \nu_{min}^{k} \geq c$ for all $k \leq r - 1$.	Bounding the quadratic form with the largest eigenvalue of the associated matrix implies
	\begin{equation}\label{eq:degen_2}
	\lim_{D_{11} \to \infty} \Big( \Big( \bd - \sum_{0 \leq k < r} x_k \bP^\trans  \bD_{k + 1} \bP  \Big)^{-1} (\bP^\trans\vh), (\bP^\trans\vh) \, \Big) \leq \big(\nu_{min}^0 \big)^{-1}  \|\vh\|^2 \leq c^{-1} \|\vh\|^2 < \infty.
	\end{equation}
	
	The logarithm terms in \eqref{eq:tracebound_transformed} can be bounded by the minimum eigenvalues in a similar manner. It suffices to show an arbitrary term in the sum is bounded, that is,
	\begin{equation}\label{eq:degen_3}
	\lim_{D_{11} \to \infty} \frac{1}{x_k}  \log \frac{|\bd - \sum_{k + 1 \leq \ell < r}x_\ell \bP^\trans \bD_{\ell + 1} \bP|}{|\bd - \sum_{k \leq \ell < r}x_{\ell} \bP^\trans \bD_{\ell + 1} \bP |} < \infty.
	\end{equation}
	If we define the matrices $\bA_k := \bd - \sum_{k \leq \ell < r}x_\ell \bP^\trans \bD_{\ell + 1} \bP$ and $\bB_k := x_{k} \bP^\trans \bD_{k + 1} \bP $, then
	\begin{equation*}
	\log \frac{|\bd - \sum_{k + 1 \leq \ell < r}x_\ell \bP^\trans \bD_{\ell + 1} \bP|}{|\bd - \sum_{k \leq \ell < r}x_{\ell} \bP^\trans \bD_{\ell + 1} \bP |} = \log \frac{|\bA_k  + \bB_k |}{|\bA_k|} = \log | \bA_k^{-1} (\bA_k  + \bB_k )| = \log | \bI + \bA_k^{-1} \bB_k |.
	\end{equation*}
	Bounding this with the largest eigenvalue, we see
	\begin{equation*}
	\log | \bI + \bA_k^{-1} \bB_k | \leq n  \log \lambda_{max}(\bI + \bA_k^{-1} \bB_k).
	\end{equation*}
	Using submultiplicativity of the spectral norm and the lower bound on the smallest eigenvalue of $\bA_k$ in \eqref{eq:smalleig}, we have
	\begin{equation*}
	\lambda_{max}(\bI + \bA_k^{-1} \bB_k) = 1 + \lambda_{max} (\bA_k^{-1} \bB_k) \leq 1 + \lambda_{max} (\bA_k^{-1}) \lambda_{max}(\bB_k) \leq 1 + c^{-1} \lambda_{max} (\bB_k) < \infty,
	\end{equation*}
	giving the required bound in \eqref{eq:degen_3}.
	
	Therefore, for a particular $\bP$, we can construct a sequence of diagonal matrices $\bd$ with arbitrary large first diagonal element such that $\frac{1}{2} \tr( \bP \bd \bP^\trans \bQ ) + F(\bP \bd \bP^\trans)$ is unbounded. In particular, we have
	\begin{equation*}
	\inf_{\bL} \Big( \frac{1}{2} \tr(\bL \bQ ) + F(\bL) \Big) = -\infty.
	\end{equation*}
	\textit{Case $\mathrm{(2)}$:} Consider the case when $\bQ$ is positive definite. We will prove that \eqref{eq:minimizefunction} attains a minimum at some point $\bL^* \in \cL$. By H\"older's inequality, $F(\bL)$ is a convex function of $\bL$, so any local minimizer is also a global minimizer. We will prove that the minimizer is attained in a compact subset of $\Gamma_n$ under the spectral norm on symmetric matrices $\| \bL \|_2 = \lambda_{max} (\bL)$.
	
	Because $\bQ$ is positive definite, the diagonal elements of $\bP^\trans \bQ \bP$ is positive and uniformly bounded away from $0$ for all orthogonal matrices $\bP$. That is, the first term in \eqref{eq:tracebound_transformed} can be bounded below by
	\begin{equation*}
	\tr\bigl( \bd \bP^\trans \bQ \bP \bigr) - \log |\bd| = \sum_{j \leq n} \bigg( \bd_{jj} (\bP^\trans \bQ \bP)_{jj} -  \log| \bd_{jj}| \bigg) \geq \sum_{j \leq n} \bigg( \bd_{jj} \lambda_{min}(\bQ) -  \log| \bd_{jj}| \bigg)  \label{eq:trace_firstterm}
	\end{equation*}
	which clearly diverges to $\infty$ if any diagonal element $\bd_{jj} \to \infty$. The remaining terms in \eqref{eq:tracebound_transformed} are non-negative, so $\frac{1}{2} \tr(\bL \bQ ) + F(\bL) \to \infty$ if $\| \bL \|_2 \to \infty$. Since $\bL - \bL_0 \geq 0$, we also have $\frac{1}{2} \tr(\bL \bQ ) + F(\bL) \to \infty$ if $\|\bL_0\|_2 \to \infty$ by monotonicity. 
	
	By definition \eqref{eq:lambda_def}, we have $\bL_{k + 1} = \bL_k + x_k \bD_{k + 1}$. By submultiplicativity of $\| \cdot \|_2$,
	\begin{equation*}
	\| \bD_{k + 1} \|_2 =  \|\bL_k \bL_k^{-1} \bD_{k + 1}  \|_2 \leq \| \bL_k\|_2 \| \bL_k^{-1} \bD_{k + 1}  \|_2,
	\end{equation*}
	so the last term in \eqref{eq:minimizefunction} can be bounded below by 
	\begin{align*}
	\log \frac{|\bL_{k + 1}|}{ |\bL_k| } = \log | \bI +  x_k \bL_k^{-1} \bD_{k + 1} | \geq \log (1 + x_k \| \bL_k^{-1} \bD_{k + 1} \|_2 ) \geq  \log(1 + x_k \| \bL_k \|_2^{-1} \| \bD_{k + 1} \|_2 ).
	\end{align*}
	Let $k^*$ be the smallest index such that $\bD_{k^* + 1} \neq \bm{0}$, then it is clear the above term diverges as $\| \bL_{k^*}\|_2 \to 0$. Since the sequence \eqref{eq:rpc2} is monotone, we have $\| \bL_{0}\|_2 = \| \bL_{k^*}\|_2$. The rest of the terms in \eqref{eq:minimizefunction} are bounded or positive for fixed $\pi$, so  $\frac{1}{2} \tr(\bL \bQ ) + F(\bL) \to \infty$ if $\| \bL_0\| \to 0$.
	
	We have shown, $\frac{1}{2} \tr(\bL \bQ ) + F(\bL)$ is unbounded if $\| \bL_0 \|_2 \to 0$ or $\| \bL_0 \|_2 \to \infty$. Therefore, there exists a $0 < c < C < \infty$ such that the minimizer is attained in the compact set
	\begin{equation*}
	\cL^* = \{ \bL \in \Gamma_n \mmm c \leq \| \bL_0 \|_2 \leq C \} \subset \cL.
	\end{equation*}
	By the extreme value theorem, there exists a $\bL_* \in \cL$ such that $\frac{1}{2} \tr(\bL \bQ ) + F(\bL)$ attains its minimum at $\bL_*$. Furthermore, our function is convex, so the minimizer $\bL_*$  is unique and satisfies the critical point condition, 
	\begin{equation*}
	\frac{\partial}{\partial t} \Big( \frac{1}{2} \tr\bigl( (\bL_* + t\bB) \bQ \bigr) + F(\bL_* + t\bB) \Big)\Big|_{t = 0} = 0,
	\end{equation*}
	for all symmetric matrices $\bB$. 
\end{proof}

Lemma \ref{lem:sharpupbd} is trivially satisfied when the infimum is $-\infty$, so we focus on the non-degenerate case moving forward. We now prove asymptotic sharpness of the upper bound using a standard large deviations calculation to decouple the constraints.

\begin{lem}\label{lem:decouple} For any $\delta^* > 0$ and positive definite $\bQ$,
	\begin{align}
	\liminf_{N \to \infty} \Phi_{\tilde{\Omega}_N^{\delta^*}}(\bL) + \frac{1}{2} \tr(\bL \bQ) \geq \inf_{\bL} \Big( \frac{1}{2} \tr(\bL \bQ) + F(\bL)  \Big) \label{eq:lowbd_rpc}.
	\end{align}
\end{lem}

\begin{proof}
	
	Consider the partition
	\begin{equation*}
	(\R^N)^n = \tilde{\Omega}_N^{\delta^*} \cup \Big( \bigcup_{j,j' \leq n} V_{j,j'}^+ \Big) \cup \Big( \bigcup_{j,j' \leq n} V_{j,j'}^- \Big)
	\end{equation*}
	where
	\begin{align}
	V_{j,j'}^+ &= \Big\{ \bvo \mmm R^{j,j'}(\bvo, \bvo) \geq Q^{j,j'} + \delta^* \Big\}, \label{eq:V+}
	\\V_{j,j'}^- &= \Big\{ \bvo \mmm R^{j,j'}(\bvo, \bvo) \leq Q^{j,j'} - \delta^* \Big\} \label{eq:V-}.
	\end{align}
	For $\bL_*$ that satisfies \eqref{eq:minimizer}, by considering values near this critical point, we will show there exists a constant $c > 0$ such that for all half-spaces $V$ in \eqref{eq:V+} or \eqref{eq:V-}
	\begin{equation}\label{eq:phiineq}
	\Phi_{V}(\bL_*) \leq F(\bL_*) - c.
	\end{equation}
	where $\Phi_{V}$ was defined in \eqref{eq:integralsubset}. We only show this for $V = V_{j,j'}^-$ for $j \neq j'$. The proof for the other cases are similar. For all $t \geq 0$ and $\bvo \in V^-_{j,j'}$, 
	\begin{equation*}
	tR^{j,j'}(\bvo, \bvo) \leq t(Q^{j,j'} - \delta^*).
	\end{equation*}
	Let $\bB$ be a matrix such that $B^{j,j'} = B^{j',j} = 1$ and is zero everywhere else. Adding and subtracting $\frac{1}{2}tNR^{j,j'}(\bvo, \bvo)$ and $\frac{1}{2}tNR^{j',j}(\bvo, \bvo)$ in the exponent, by symmetry of $\bQ$, we have
	\begin{align}
	\Phi_{ V }(\bL_*) &\leq t (Q^{j,j'} - \delta^*) + \Phi_{V} (\bL_* + t \bB) \notag
	\\&\leq t (Q^{j,j'} - \delta^*) + F (\bL_* + t\bB) \notag
	\\&= -t\delta^* - \frac{1}{2} \tr(\bL_* \bQ) + \frac{1}{2} \tr((\bL_* + t \bB) \bQ ) + F(\bL_* + t\bB) =:U(t) \label{eq:upbdt}.
	\end{align}
	Since $U(0) = F(\bL_*)$, the critical point condition \eqref{eq:minimizer} implies $U'(0) = -\delta^*$. In particular, there is a $t^*$ such that $U(t^*) < U(0)$. Since \eqref{eq:upbdt} holds for all $t > 0$, there is a $c$ such that
	\begin{equation*}
	\Phi_{V}(\bL_*) \leq U(t^*) \leq U(0) - c = F(\bL_*) - c.
	\end{equation*}

	Recall the sets \eqref{eq:V+}, \eqref{eq:V-}, and \eqref{eq:V0} form a partition of $(\R^N)^n$. A consequence of the recursion in the Ruelle probability cascades (see equation (118) in the proof of \cite[Lemma 7]{PVS}) implies
	\begin{equation*}
	F(\bL_*) \leq \frac{\log (2n^2 + 1)}{N x_0} + \max\Bigl( \max_{V} \Phi_{V}(\bL_*), \Phi_{\tilde{\Omega}_N^{\delta^*}}(\bL_*) \Bigr)
	\end{equation*}
	where the maximum over $V$ is over the halfspaces of the form \eqref{eq:V+}, \eqref{eq:V-}. Our bounds in \eqref{eq:phiineq} ensures we cannot have
	\begin{equation*}
	F(\bL_*) \leq \frac{\log (2n^2 + 1)}{N x_0} + \max_{V} \Phi_{V}(\bL_*)
	\end{equation*}
	for $N$ sufficiently large. Therefore, we must have
	\begin{equation*}
	F (\bL_*) \leq \frac{\log (2n^2 + 1)}{N x_0} +  \Phi_{\tilde{\Omega}_N^{\delta^*}} (\bL_*).
	\end{equation*}
	Taking $N \to \infty$ completes the proof.
\end{proof}

The proof of Lemma \ref{lem:sharpupbd} follows by applying Lemma \ref{lem:decouple} to \eqref{eq:boundapprox},
\begin{equation*}
\lim_{\epsilon \to 0} \liminf_{N \to \infty} f_N^1(\pi) \geq  \inf_{\bL} \frac{1}{2} \bigg( \tr(\bL \bQ) - n - \log |\bL|  + (\bL^{-1}_0 \vh,\vh) +  \sum_{0 \leq k \leq r-1} \frac{1}{x_k}  \log \frac{|\bL_{k + 1}|}{|\bL_k|} \bigg).
\end{equation*}

We have shown that Theorem \ref{theo:MAIN} is trivially satisfied for degenerate constraint $\bQ$ just from examining the upper bound. The case for positive definite constraint $\bQ$ is much harder and will require some preliminary work before attempting the cavity computations. We begin by introducing a variant of the Aizenman--Sims--Starr scheme.

\section{The Aizenman--Sims--Starr Scheme}\label{sec:ASS}

Before we can complete the cavity computations to prove the lower bound, we first prove an analogue of the Aizenman--Sims--Starr scheme \cite{AS2} for spherical spin glass models with vector spins. The extension to this model is non-trivial because the uniform measure on the sphere is not a product measure, so the usual proof of the scheme fails. 

This section follows the proof of the Aizeman--Sims--Starr scheme adapted for spherical models in \cite{CASS}. The main difference is the Aizenman--Sims--Starr scheme (see Lemma~\ref{lem:ASS_Pert}) will be with respect to a Gaussian reference measure as opposed to the surface measure in \cite{CASS}. This form was chosen for convenience, because it matches the form of the functional \eqref{eq:upbd4}. 

To simplify notation, we first prove an analogue of the Aizenman--Sims--Starr scheme with no external field. We will explain how to reintroduce the external field at the end of Section~\ref{sec:modpert}. Consider the partition function with $\veh = \vec{0}$ for a system of size $N$,
\begin{equation}\label{SHOWLABEL}
Z_N(\bQ,\epsilon) = \int_{Q^\epsilon_{N}} \exp \big( H_N(\bvs) \big) \, d \lambda_N^n(\bvs)
\end{equation}
and the corresponding partition function for a system of size $M + N$,
\begin{equation}\label{SHOWLABEL}
Z_{M + N}(\bQ,\epsilon) = \int_{Q^\epsilon_{M + N}} \exp \big( H_{N+M}(\bvp) \big) \, d \lambda_{M + N}^n(\bvp).
\end{equation}
We denote spin configurations from the system of size $M + N$ with $\bvr = (\bvs, \bvo) \in S^{n}_{M + N}$ where $\bvs \in \R^N$ denotes the bulk coordinates and $\bvo \in \R^M$ denotes the cavity coordinates. 

We proceed like the traditional Aizenman--Sims--Starr scheme and split the Hamiltonian into the cavity fields \cite[Section 3.5]{PBook}
\begin{gather}
H_{M + N} (\bvs, \bvo) \eqdist \sum_{j \leq n} {H}^j_{M,N}(\bvs) + \sum_{i \leq M} \sum_{j \leq n} \vo_i(j) Z_i^j(\bvs) + r(\bvp),\label{eq:normal}\\
H_{N} (\bvs, \bvo) \eqdist \sum_{j \leq n} {H}^j_{M,N}(\bvs) + \sqrt{M} \sum_{j \leq n} Y^j(\bvs).\label{eq:cavity}
\end{gather}
Here, $H_{M,N}$ is defined like $H_N$ but with normalization $(M + N)^{-(p-1)/2}$. The covariance of this Hamiltonian is given by
\begin{equation}
\E {H}^{j}_{M,N}(\bvs^\ell) {H}^{j'}_{M,N}(\bvs^{\ell'}) = (M + N) \xi_{j,j'} \Bigl( \frac{N}{M + N}  R_{\ell, \ell'}^{j,j'} \Bigr).\label{eq:fieldH}
\end{equation}
The cavity fields $Z(\bvs)$ and $Y(\bvs)$ in \eqref{eq:normal} and \eqref{eq:cavity} are centered Gaussian processes with covariances:
\begin{align}
\E Z_i^j(\bvs^\ell) Z_{i'}^{j'}(\bvs^{\ell'}) &= \delta_{i,i'}\xi'_{j,j'} \bigl( R_{\ell, \ell'}^{j,j'} \bigr) + \sO \Big( \frac{M}{N} \Big),\label{eq:fieldZ}\\
\E Y^j(\bvs^\ell) Y^{j'}(\bvs^{\ell'}) &= \theta_{j,j'} \bigl( R_{\ell, \ell'}^{j,j'} \bigr) + \sO \Big( \frac{M}{N} \Big),\label{eq:fieldY}
\end{align}
and the remainder term $r(\bvp)$ has covariance,
\begin{equation}
\E r(\bvp^{\,\ell}) r(\bvp^{\,\ell'}) = \sO \Big( \frac{M^2}{M + N} \Big)\label{eq:fieldr}.
\end{equation}
We will prove that we can replace the cavity fields $Z(\bvs)$ and $Y(\bvs)$ with centered Gaussian fields $z_i(\bvs)$ and $y(\bvs)$ taking values in $\R^n$ indexed by $\bvs \in S_N^n$, with covariances
\begin{align}
\E z_i^j(\bvs^\ell) z_{i'}^{j'}(\bvs^{\ell'}) &= \delta_{i,i'}\xi'_{j,j'} \bigl( R_{\ell, \ell'}^{j,j'} \bigr),\label{eq:fieldz}\\
\E y^j(\bvs^\ell) y^{j'}(\bvs^{\ell'}) &= \theta_{j,j'} \bigl( R_{\ell, \ell'}^{j,j'} \bigr)\label{eq:fieldy}.
\end{align}

Let $\langle \cdot \rangle_{M,N}$ be the average with respect to the Gibbs measure,
\begin{equation}\label{SHOWLABEL}
G_{M,N}(d\bvs) = \frac{\exp H_{M,N}(\bvs) \, d \lambda_N^n(\bvs)}{Z_{M,N}(\bQ, \epsilon) },
\end{equation}
on $Q^\epsilon_N$, with normalization
\begin{equation}\label{SHOWLABEL}
Z_{M,N}(\bQ,\epsilon) = \int_{Q^\epsilon_{N}} \exp \big( H_{M,N}(\bvs) \big) \, d \lambda_{N}^n(\bvs).
\end{equation}

We start as usual with the inequality
\begin{align}
\liminf_{N \to \infty} \frac{1}{N} \E \log Z_N(\bQ,\epsilon) &\geq \frac{1}{M} \liminf_{N \to \infty} \big( \E \log Z_{M + N}(\bQ,\epsilon)  - \E \log Z_{N}(\bQ,\epsilon) \big) \notag
\\&= \frac{1}{M} \liminf_{N \to \infty} \bigg( \E \log \frac{Z_{M + N}(\bQ,\epsilon)}{Z_{M,N}(\bQ,\epsilon) }  - \E \log \frac{Z_{N}(\bQ,\epsilon)}{Z_{M,N}(\bQ,\epsilon)} \bigg) \label{eq:ass_step_2}.
\end{align}
The surface measure $\lambda_{M + N}$ appearing in $Z_{M + N}$ is not a product measure, so the standard proof of the Aizenman--Sims--Starr does not apply in this setting. We will prove the Aizenman--Sims--Starr scheme for the spherical spin glass model with vector spins. Let the $\delta$ shell around $Q_M^\epsilon$ be denoted by
\begin{equation}\label{eq:deltashell}
\Omega_M^{\epsilon,\delta} = \big\{  \bvo = (s_j \bt(j))_{j \leq n} \in (\R^M)^n \mmm \bvt \in Q_M^{\epsilon}, \, s_j \in [\sqrt{1 - \delta}, \sqrt{1 + \delta}] \text{ for all } j \leq n \big\},
\end{equation}
where $s_j \in \R^+$ and $\tau(j) \in S_M$ are the radial and angular components of the polar form of $\bo(j)$ (see after \eqref{eq:upbd_Y_recur} for the formulas). We will prove the following lower bound of \eqref{eq:ass_step_2}.

\begin{lem}\label{lem:ass_unpert}
	Let $\nu_M$ be the standard normal distribution on $\R^M$. For any $\epsilon > 0$ and $M \geq 1$, there exists a $\delta \in (0,\epsilon)$ such that \eqref{eq:ass_step_2} is bounded below by
	\begin{align}
	&\frac{1}{M} \liminf_{N \to \infty} \bigg( \E \log \Big \langle \int_{\Omega_M^{\epsilon / 2, \delta}} \exp \Big( \sum_{i \leq M} \sum_{j \leq n} \vo_i(j) z_i^j(\bvs) \Big) \, d \nu_{M}^n(\bvo) \Big \rangle_{M,N} - \E \log \Big\langle  \exp \sqrt{M} y(\bvs) \Big\rangle_{M,N} \Big) \bigg) \notag\\
	&\qquad - L \delta. \label{eq:ass_error}
	\end{align}
\end{lem}
The main difference between the bound \eqref{eq:ass_error} and the traditional Aizenman--Sims--Starr representation is the Gaussian reference measure appearing in the first cavity field. This measure appears as a consequence of the Poincar\'e limit, which states that the standard Gaussian measure in $\R^M$ is the limiting distribution of projected uniform distributions on $S_{N + M}$ as $N$ tends to infinity. 

\subsection{Poincar\'e limit} We first explain a method to asymptotically decouple $\lambda_{M + N}$ into an approximate product measure over the spheres $S_N \times S_M$. The distribution of the projection of $S_{N + M}$ onto $\R^M$ under $\lambda_{N + M}$ converges weakly to the Gaussian distribution $\nu_M$ on $\R^M$ in the Poincar\'e limit \cite{poincare1912calcul}. In particular, the distribution of the cavity coordinates under the normalized surface measure will be approximately Gaussian for large $N$. For large $M$, $\nu_M$ will concentrate around $S_M$. We first introduce some notation and state this result in one dimension.

For $K \geq 1$, we denote the unit sphere in $\R^K$ with $S_K^1$ and $|S_K^1|$ its surface area. Let
\begin{align}\label{SHOWLABEL}
A_{M,N} &= \prod_{j = 1}^M \Big[ - \sqrt{M + N + 1 - j}, \sqrt{M + N + 1 - j} \,\Big],
\end{align}
be a subset of $\R^M$ representing the domain of the cavity coordinates. We define the density on $\R^M$,
\begin{align*}
d \nu_{M,N} (\bm{x})= f_{M,N}(\bm{x}) \, d \bm{x},
\end{align*}
where
\begin{equation}\label{eq:densityFMN}
f_{M,N}(\bm{x}) = b_{M,N} \prod_{j=1}^{M} \bigg(1 - \frac{x_j^2}{M + N + 1 - j}\bigg)^{\frac{M + N - j - 2}{2}},
\end{equation}
with normalizing coefficient
\begin{equation*}
b_{M,N} = \prod_{j = 1}^M \frac{|S_{M + N - j}^1| }{ |S^1_{M + N + 1 - j}| \sqrt{M + N + 1 - j}}.
\end{equation*}
The pointwise limit of \eqref{eq:densityFMN} converges to the standard normal distribution on $\R^M$
\begin{align*}
d \nu_{M} (\bm{x})= f_{M}(\bm{x}) \, d \bm{x},
\end{align*}
where
\begin{equation*}
f_M(\bm{x}) := \lim_{N \to \infty} f_{M,N} (\bm{x}) = \bigg( \frac{1}{2\pi} \bigg)^{M / 2} \exp \biggl( - \frac{\| \bm{x} \|^2 }{2} \biggr).
\end{equation*}
Lastly, we define the coefficients
\begin{equation}\label{eq:coeffs_a}
a_1 = 1, \quad a_\ell(\bm{x}) = \prod_{j=1}^{\ell -1} \sqrt{1 + \frac{1 - x_j^2}{M + N - j}} \text{ for $1 < \ell \leq M + 1$}
\end{equation}
and the corresponding map for $\psi: S_N \times A_{M,N} \to S_{M + N}$ given by
\begin{equation*}
\psi(\bm{\sigma}, \bm{\omega}) = \big(\sigma_1 a_{M + 1}(\bo) ,\dots,\sigma_Na_{M + 1}(\bo), \omega_1 a_1(\bo), \dots, \omega_M a_M(\bo) \big)
\end{equation*}
for $\bs \in S_N$ and $\bo \in A_{M,N}$. The surface measure on $S_{M + N}$ can be decoupled as follows:
\begin{lem}\label{lem:poincare_1d} {\cite[Lemma 3]{CASS}}
	Suppose $g$ is a nonnegative function defined on $S_{M + N}$.  Then for $\bm{\rho} = (\bm{\sigma}, \bm{\omega}) \in S_{M + N}$ we have
	\begin{equation*}
	\int_{S_{M + N}} g(\bm{\rho}) \, d \lambda_{M + N} (\bm{\rho}) = \int_{A_{M,N}} \int_{S_N}    g( \psi(\bm{\sigma}, \bm{\omega}) ) \,d\lambda_N(\bm{\sigma})  d \nu_{M,N} (\bo) .
	\end{equation*}
\end{lem}

We will need a multidimensional version of this argument. To simplify notation, for $n$ copies of $S_{N + M}$, we define $a_{\ell}^j := a_{\ell} (\bo(j))$ keeping the dependence on the cavity coordinate $\bo(j)$ implicit. Similarly, we define
\begin{equation*}
\Psi(\bvs, \bvo) = \Big( \psi\big(\bs(j), \bo(j)\big) \Big)_{j \leq n} = \Big( a_{M + 1}^j \vs_1(j),\dots,a_{M + 1}^j \vs_N(j), a_1^j \vo_1(j) , \dots, a_M^j \vo_M(j)  \Big)_{j \leq n},
\end{equation*}
to represent the transformation applied coordinate-wise. The following result explains how the surface measure on $S_{M + N}$ decouples asymptotically. 
\begin{cor}\label{cor:poincare_lemma}
	Suppose $g$ is a nonnegative function defined on $S^n_{M + N}$.  Then for $\bvp = (\bvs, \bvo) \in S^n_{M + N}$ we have
	\begin{equation}\label{eq:poincare_lemma}
	\int_{S^n_{M + N}} g(\bvp) \, d \lambda^n_{M + N} (\bvp) = \int_{S^n_N} \int_{A^n_{M,N}}  g\big( \Psi( \bvs, \bvo ) \big) \,d \nu^n_{M,N} (\bvo) d\lambda^n_N(\bvs)   .
	\end{equation}
\end{cor}  
\begin{proof} We apply Lemma \ref{lem:poincare_1d} to each coordinate $j \leq n$. The region of integration is a product set and we are integrating a non-negative function, so we are freely able to rearrange the order of integration by Fubini's Theorem. 
\end{proof}
We apply Corollary \ref{cor:poincare_lemma} to lower bound $Z_{M + N}(\bQ, \epsilon)$ with an integral over the product set of bulk and cavity coordinates. To simplify notation, we denote the transformed coordinates with 
\begin{equation}\label{eq:poincarecoordsrho}
\ubvp = (\ubvs, \ubvo) := \Psi( \bvs, \bvo)
\end{equation}
where 
\begin{equation}\label{eq:poincarecoords}
\ubvs = \big(a^j_{M + 1}\vs_1(j), \dots, a^{j}_{M + 1}\vs_N(j) \big)_{j \leq n} \text{ and } \ubvo = \big(a^j_1 \vo_1(j), \dots, a^j_M \vo_M(j) \big)_{j \leq n}
\end{equation}
are the respective transformed bulk and cavity coordinates. For an arbitrary non-negative function $g$ on $S^n_{M + N}$, \eqref{eq:poincare_lemma} implies
\begin{align}
\int_{Q_{M + N}^\epsilon} g(\bvp) \, d \lambda_{M + N}^n(\bvp) &= \int_{S^n_{M + N}} \1_{Q_{M + N}^\epsilon} (\bvp) g(\bvp) \, d\lambda_{M + N}^n (\bvp) \notag\\
&= \int_{S^n_N} \int_{A^n_{M,N}}  \1_{Q_{M + N}^\epsilon}( \ubvs, \ubvo ) g( \ubvs, \ubvo ) \,  d \nu^n_{M,N} (\bvo) d\lambda^n_N(\bvs)\label{eq:poincare_indicator}.
\end{align}
We first split the integral over $Q^\epsilon_{M + N}$ into a product set over $Q_N^\epsilon \times \Omega_M^{\epsilon/2,\delta}$ for suitably chosen $\delta$.
\begin{lem}\label{lem:split_cavity}
	For any $\epsilon > 0$ and $N$ sufficiently large, there exists a $\delta \in (0,\epsilon)$ such that
	\begin{equation}\label{eq:splitconstraint}
	\1_{Q_{M + N}^\epsilon} ( \ubvs,\ubvo) \geq \1_{Q_N^\epsilon} (\bvs) \1_{\Omega_M^{\epsilon/2,\delta}} (\bvo).
	\end{equation}
\end{lem}
\begin{proof}
	We will find conditions on $\delta$ such that \eqref{eq:splitconstraint} holds. Let $\delta > 0$ and take $\bvs \in Q_N^\epsilon$ and $\bvo \in \Omega_M^{\epsilon/2,\delta}$. The overlaps of the transformed coordinates $\ubvp$ defined in \eqref{eq:poincarecoordsrho} satisfy
	\begin{equation}\label{eq:overlaprhosplit}
	\bR( \ubvp, \ubvp ) = \frac{N }{M + N}\bR(\ubvs, \ubvs) + \frac{M }{M + N}\bR(\ubvo, \ubvo).
	\end{equation}
	The set $\Omega_M^{\epsilon/2,\delta}$ is bounded, so the corresponding transformed overlaps $\bR(\ubvs, \ubvs)$ and $\bR(\ubvo, \ubvo)$, can be approximated by the standard overlaps $\bR(\bvs,\bvs)$ and $\bR(\bvt, \bvt)$ of configurations $\bvs$ and $\bvt$ on the spheres $S_N$ and $S_M$ respectively. 
	
	Firstly, since $\|\bvo(j)\|^2 < M (1 + \delta)$ for all $\bvo \in \Omega_M^{\epsilon/2,\delta}$ and $j \leq n$, we have the relation
	\begin{equation}\label{eq:limitcoeffN}
	\lim_{N \to \infty} \big(N - a_{M + 1}^j a_{M + 1}^{j'} N \big) = \frac{\|\bo(j)\|^2 + \|\bo(j')\|^2}{2} - M \leq M\delta.
	\end{equation}
	Since $R^{j,j'}(\ubvs,\ubvs) = a_{M + 1}^j a_{M + 1}^{j'} R^{j,j'}(\bvs,\bvs)$, for all $N$ sufficiently large
	\begin{equation}\label{eq:transformedoverlapN}
	\Big \| \frac{N}{N + M} \bR(\ubvs, \ubvs) - \frac{N}{N + M} \bR(\bvs, \bvs) \Big\|_\infty  \leq \frac{ L M \delta}{N + M}.
	\end{equation}
	
	Secondly, $\Omega_M^{\epsilon/2,\delta}$ is a compact set so $\lim_{N \to \infty }a_\ell(\bvo) = 1$ uniformly for all $1 < \ell \leq M$. Therefore,
	\begin{equation*}
	\lim_{N \to \infty} \| \bR(\ubvo, \ubvo) - \bR(\bvo, \bvo) \|_\infty = 0,
	\end{equation*}
	uniformly on $\Omega_M^{\epsilon/2,\delta}$. Likewise, on $\Omega_M^{\epsilon/2,\delta}$, $\| \bvo(j) - \bvt(j) \|^2 = \| s_j\bvt(j) - \bvt(j) \|^2 \leq \delta M$ for all $j \leq n$, so
	\begin{equation*}
	\| \bR(\bvo, \bvo) - \bR(\bvt, \bvt) \|_\infty \leq \delta.
	\end{equation*}
	Therefore, the triangle inequality implies for all $\bvo \in \Omega_M^{\epsilon/2,\delta}$,
	\begin{equation}\label{eq:transformedoverlapM}
	\Big \| \frac{M}{N + M} \bR(\ubvo, \ubvo) - \frac{M}{N + M} \bR(\bvt, \bvt) \Big\|_\infty \leq \frac{LM \delta}{N + M}.
	\end{equation}
	
	For any $\bvs \in Q_N^\epsilon$ and $\bvo \in \Omega_M^{\epsilon/2,\delta}$, \eqref{eq:transformedoverlapN} and \eqref{eq:transformedoverlapM} imply for $N$ sufficiently large
	\begin{align*}
	\big\|\bR(\ubvp,\ubvp) - \bQ \big\|_\infty &=  \Big\| \frac{N }{M + N} \big(\bR(\ubvs, \ubvs) - \bQ \big) + \frac{M }{M + N} \big(\bR(\ubvo, \ubvo) - \bQ \big) \Big\|_\infty
	\\&\leq \Big\| \frac{N }{M + N} \big(\bR(\bvs, \bvs) - \bQ \big) + \frac{M }{M + N} \big(\bR(\bvt, \bvt) - \bQ \big) \Big\|_\infty + \frac{L M \delta}{N + M}
	\\&\leq \frac{N\epsilon}{N + M} + \frac{M\epsilon}{2(N + M)}  + \frac{L M \delta}{M + N}
	\\&= \epsilon + \frac{M}{N + M} \bigg(  L \delta -\frac{\epsilon}{2} \bigg).
	\end{align*}
	Choosing $\delta \leq \frac{\epsilon}{2 L}$, we have
	\begin{equation*}
	\big\|\bR(\ubvp,\ubvp) - \bQ \big\|_\infty \leq \epsilon.
	\end{equation*}
	In particular, this means
	\begin{equation*}
	\Big\{ (\ubvs,\ubvo) \in Q_{N + M}^\epsilon \Big\} \supseteq \Big\{ \bvs \in Q^\epsilon_N  \Big\} \times \Big\{ \bvo \in  \Omega_M^{\epsilon/2,\delta} \Big\},
	\end{equation*}
	completing the proof.
\end{proof}

Applying Lemma \ref{lem:split_cavity} to \eqref{eq:poincare_indicator} and taking $g(\bvp) = \exp H_{M + N}(\bvp)$, we have for $N$ sufficiently large
\begin{equation}\label{eq:step1_sub}
\E \log Z_{M + N}(\bQ,\epsilon) \geq \E \log \int_{Q_N^{\epsilon}} \int_{\Omega_M^{\epsilon/2,\delta}} \exp \Big(\sum_{j \leq n} H^j_{M + N} \big( \ubvs,\ubvo \big)  \Big) \, d\nu_{M,N}^n(\bvo) d \lambda^n_N(\bvs).
\end{equation}
Consequently, we are able to decouple the surface measure, which resolves the first major obstacle in the proof of the Aizenman--Sims--Starr representation.

\subsection{Proof of Lemma~\ref{lem:ass_unpert}} Using \eqref{eq:step1_sub}, we can derive a lower bound for the first term in \eqref{eq:ass_step_2}.

\begin{lem}\label{lem:ASSfirstterm}
	For every $\epsilon > 0$, there exists a $\delta \in (0,\epsilon)$ such that
	\begin{align}
	\liminf_{N \to \infty} \E \log \frac{Z_{M + N}(\bQ,\epsilon)}{Z_{M,N}(\bQ,\epsilon)} &\geq \liminf_{N \to \infty} \E \log \Big \langle \int_{\Omega_M^{\epsilon/2,\delta}} \exp \Big( \sum_{i \leq M} \sum_{j \leq n} \vo_i(j) z_i^j(\bvs) \Big) \, d \nu_{M}^n(\bvo) \Big \rangle_{M,N}  \notag
	\\&\qquad\quad\quad  - L M \delta \label{eq:assfirstterm}.
	\end{align} 
\end{lem}

\begin{proof}
	We start by splitting the left hand side of \eqref{eq:assfirstterm} in three parts
	\begin{align}\label{eq:ass_split}
	\!\!\!\!\!
	\E \log \dfrac{Z_{M + N}(\bQ,\epsilon)}{\int_{\Omega_M^{\epsilon/2,\delta}} J_{M,N} \, d{\nu}_{M,N}^n(\bvo) } + \E \log \dfrac{\int_{\Omega_M^{\epsilon/2,\delta}} J_{M,N} \, d{\nu}_{M,N}^n(\bvo) }{\int_{\Omega_M^{\epsilon/2,\delta}} J_{M,N} \, d{\nu}_{M}^n(\bvo) } + \E \log \dfrac{\int_{\Omega_M^{\epsilon/2,\delta}} J_{M,N} \, d{\nu}_{M}^n(\bvo) }{Z_{M,N}(\bQ,\epsilon)},
	\end{align}
	where
	\begin{align*}
	J_{M,N} &:= J_{M,N}(\bvo) = \int_{Q_N^{\epsilon}} \exp \Big(\sum_{j \leq n} H^j_{M,N} (\bvs) + \sum_{j \leq n} \sum_{i \leq M} \vo_i(j) z_i^j(\bvs) \Big) \, d \lambda^n_N(\bvs).
	\end{align*}
	We bound each of the terms in \eqref{eq:ass_split} separately. 
	\\\\\textit{Step 1:} We show the first term in \eqref{eq:ass_split} satisfies
	\begin{align}\label{eq:ASS_Step1}
	\liminf_{N \to \infty} \E \log \dfrac{Z_{M + N}(\bQ,\epsilon)}{\int_{\Omega_M^{\epsilon/2,\delta}} J_{M,N} \, d{\nu}_{M,N}^n(\bvo) } \geq - LM\delta . 
	\end{align}
	Recall the lower bound \eqref{eq:step1_sub}
	\begin{equation}\label{eq:lowbd_repeat}
	\E \log Z_{M + N}(\bQ,\epsilon) \geq \E \log \int_{Q_N^{\epsilon}} \int_{\Omega_M^{\epsilon/2,\delta}} \exp \Big(\sum_{j \leq n} H^j_{M + N} \big( \ubvs,\ubvo \big)  \Big) \, d\nu_{M,N}^n(\bvo) d \lambda^n_N(\bvs).
	\end{equation}
	We now use Gaussian interpolation to control the term on the right hand side. 
	
	Recall the Gaussian fields in \eqref{eq:normal} and \eqref{eq:fieldz}. We define the interpolating Hamiltonian
	\begin{equation}\label{eq:triplehamiltonian}
	H_t(\bvp) = \sum_{j \leq n}H^j_{1,t}(\bvp) + H^j_{2,t} (\bvp) + H^j_{3,t} (\bvp).
	\end{equation}
	For $\bvs \in Q_N^\epsilon$, $\bvo \in \Omega_M^{\epsilon/2,\delta}$, using the corresponding transformed coordinates \eqref{eq:poincarecoords}, the Gaussian processes in \eqref{eq:triplehamiltonian} are given by
	\begin{align}
	H^j_{1,t}(\bvp) &= \sqrt{t} H^j_{M,N} \big(\ubvs \big) + \sqrt{1 - t} \tilde{H}^j_{M,N} (\bvs) \label{eq:ham_1},
	\\H^j_{2,t}(\bvp) &= \sum_{i \leq M} \vo_i(j) \Big( \sqrt{t} a_i^j Z^j_i(\ubvs ) + \sqrt{1 - t} z^j_i(\bvs)   \Big) \label{eq:ham_2},
	\\H^j_{3,t}(\bvp) &= \sqrt{t} \, r\big(\ubvs,\ubvo \big) \label{eq:ham_3}.
	\end{align}
	The Gaussian process $\tilde{H}^j_{M,N}$ is an independent copy of $H^j_{M,N}$. Let
	\begin{equation}\label{eq:interpolatingfunction3}
	\phi(t) = \E \log \int_{Q_N^\epsilon} \int_{\Omega_M^{\epsilon/2,\delta}} \exp \Bigl(H_t (\bvs, \bvo)  \Bigr) \, d\nu_{M,N}^n(\bvo) d \lambda^n_N(\bvs)
	\end{equation}
	be the corresponding interpolating Hamiltonian. By Gaussian integration by parts,
	\begin{equation}\label{eq:cov_ass}
	\phi'(t) = \frac{1}{2} \sum_{j,j' \leq n} \sum_{k \leq 3} \E \bigg\langle \E \frac{\partial H^j_{k,t} (\bvp^{\,1})}{ \partial t} \cdot H^{j'}_{k,t}(\bvp^{\,1}) - \E \frac{\partial H^j_{k,t} (\bvp^{\,1})}{ \partial t} \cdot H^{j'}_{k,t}(\bvp^{\,2}) \bigg\rangle_t
	\end{equation}
	where $\langle \cdot \rangle_t$ is the average with respect to the Gibbs measure on $Q_N^\epsilon \times \Omega_M^{\epsilon/2,\delta}$ proportional to $\exp(H_t)$ with respect to the reference measure $\lambda^n_N \times \nu_{M,N}^n$. We now compute the covariances of the cavity fields in \eqref{eq:triplehamiltonian}. 
	
	For all $\ell \leq M + 1$ and $j \leq n$, $a^j_{\ell}(\bvo) \to 1$ uniformly on $\Omega_M^{\epsilon/2,\delta}$ by compactness. The leading terms of \eqref{eq:fieldZ} and \eqref{eq:fieldr} do not grow in $N$, so by continuity the terms \eqref{eq:ham_2} and \eqref{eq:ham_3} in \eqref{eq:cov_ass} vanish in the limit. 

	We now compute the covariances containing $H^j_{1,t}$. The covariance of \eqref{eq:fieldH} is order $N$, so it is not obvious that the differences of the covariances are small. We resolve this by using identity \eqref{eq:limitcoeffN} and applying the mean value theorem. If we let $a^j_{i}(\bo^1) := a_{i}(\bo^1(j))$, then
	\begin{align*}
	&\lim_{N \to \infty} \Big|\E \frac{\partial H^j_{1,t} (\bvp^{\,1})}{ \partial t} \cdot H^{j'}_{1,t}(\bvp^{\,2})  \Big|\\
	&= \lim_{N \to \infty} (M + N) \Big| \xi'_{j,j'} \Big( \frac{N}{M + N} a_{M + 1}^{j}(\bo^1)a_{M + 1}^{j'}(\bo^2) R_{1,2}^{j,j'} \Big) - \xi'_{j,j'} \Big( \frac{N}{M + N}R_{1,2}^{j,j'} \Big) \Big|  \\
	&\leq \lim_{N \to \infty} \|\bxi'(1)\|_1  \bigg|  N \Big( a^j_{M + 1}(\bo^1)a^{j'}_{M + 1}(\bvo^2) - 1 \Big) \bigg|\\ 
	&\leq  \|\bxi'(1)\|_1 \delta M,
	\end{align*}
	uniformly on $Q_N^\epsilon \times \Omega_M^{\epsilon/2,\delta}$ for all $j,j' \leq n$. In particular, $\lim_{N \to \infty} \sup_{0 \leq t \leq 1} |\phi'(t)| \leq L M \delta$.  The mean value theorem implies
	\begin{align*}
	\liminf_{N \to \infty} \E \log \dfrac{Z_{M + N}(\bQ,\epsilon)}{\int_{\Omega_M^{\epsilon/2,\delta}} J_{M,N} \, d{\nu}_{M,N}^n(\bvo) } &\geq \liminf_{N \to \infty} \big( \phi(1) - \phi(0) \big) \geq - \limsup_{N \to \infty} \sup_{0 \leq t \leq 1} | \phi'(t)| \geq  -  L M \delta,
	\end{align*}
	finishing the bound of \eqref{eq:ASS_Step1}.
	\\\\\textit{Step 2:} We show the second term in \eqref{eq:ass_split} satisfies
	\begin{align*}
	\liminf_{N \to \infty} \E \log \dfrac{\int_{\Omega_M^{\epsilon/2,\delta}} J_{M,N} \, d{\nu}_{M,N}^n(\bvo) }{\int_{\Omega_M^{\epsilon/2,\delta}} J_{M,N} \, d{\nu}_{M}^n(\bvo) } \geq 0.
	\end{align*}
	This proof is identical to the proof of Lemma 5 in \cite{CASS}. The key observation is $f_{M,N}$, the density of $\nu_{M,N}$, converges to $f_M(\bvo)$, the density of $\nu_M$. Since $\log(1 - x) \geq -x - x^2$ for $x < 0.5$, for $N$ sufficiently large, we have
	\begin{align*}
	\log \frac{f_{M,N}( \bm{x}) }{f_{M}(\bm{x}) } &= \log b_{M,N} + \sum_{j = 1}^M \frac{M + N - 2 - j}{2} \cdot  \log \Bigl(1 - \frac{x_j^2}{M + N + 1 - j} \Bigr) - \frac{M}{2}\log \frac{1}{2\pi} + \frac{\| \bm{x} \|^2}{2} \notag
	\\& \geq \log b_{M,N} + \sum_{j = 1}^M \bigg( \frac{M + N  + 1 - j}{2}  \bigg) \cdot \bigg(-\frac{x_j^2}{M + N + 1 - j} - \frac{x_j^4}{(M + N + 1 - j)^2}  \bigg)  
	\\& \quad - \frac{M}{2}\log \frac{1}{2\pi} + \frac{\| \bm{x} \|^2}{2} \notag \color{black} 
	\\&\geq \log \frac{b_{M,N}}{ (2 \pi)^{-M/2}} -   \frac{ \|\bm{x}\|^4 }{2N}. 
	\end{align*}
	
	To simplify notation, let $f^n_{M}(\bvo) := \prod_{j \leq n} f_{M}(\bvo(j))$ and $f^n_{M,N}(\bvo) := \prod_{j \leq n} f_{M,N}(\bvo(j))$.   Jensen's inequality implies
	\begin{align*}
	\E \log \dfrac{\int_{\Omega_M^{\epsilon/2,\delta}} J_{M,N} \, d{\nu}_{M,N}^n(\bvo) }{\int_{\Omega_M^{\epsilon/2,\delta}} J_{M,N} \, d{\nu}_{M}^n(\bvo) } &=\E \log \frac{\int_{\Omega_M^{\epsilon/2,\delta}} J_{M,N}(\bvo)f^n_{M}(\bvo) \frac{  f^n_{M,N}(\bvo)}{f^n_{M}(\bvo) } d\bvo }{\int_{\Omega_M^{\epsilon/2,\delta}} J_{M,N}(\bvo) f^n_M (\bvo) d\bvo}\\
	&\geq \E \frac{\int_{\Omega_M^{\epsilon/2,\delta}} J_{M,N}(\bvo) f^n_M (\bvo)   \log \frac{f^n_{M,N}(\bvo) }{f^n_{M}(\bvo) } d\bvo }{\int_{\Omega_M^{\epsilon/2,\delta}} J_{M,N}(\bvo) f^n_M (\bvo) \,  d\bvo}\\
	&\geq n\log \frac{b_{M,N}}{ (2 \pi)^{-M/2}} -  \frac{n (1 + \delta)^2 M^2 }{2N}.
	\end{align*}
	Since $b_{M,N} \to (2\pi)^{-M/2}$, we have 
	\begin{equation*}
	\liminf_{N \to \infty}\E \log \dfrac{\int_{\Omega_M^{\epsilon/2,\delta}} J_{M,N} \, d{\nu}_{M,N}^n(\bvo) }{\int_{\Omega_M^{\epsilon/2,\delta}} J_{M,N} \, d{\nu}_{M}^n(\bvo) } \geq \liminf_{N \to \infty} n\log \frac{b_{M,N}}{ (2 \pi)^{-M/2}} - \frac{n (1 + \delta)^2 M^2 }{2N} = 0.
	\end{equation*}	
	\color{black}
	\\\\\textit{Step 3:} By definition of $\langle \cdot \rangle_{M,N}$, the last term in \eqref{eq:ass_split} is equal to
	\begin{align*}
	\liminf_{N \to \infty} \E \log \dfrac{\int_{\Omega_M^{\epsilon/2,\delta}} J_{M,N} \, d{\nu}_{M}^n(\bvo) }{Z_{M,N}(\bQ,\epsilon)} &= \liminf_{N \to \infty} \E \log \Big \langle \int_{\Omega_M^{\epsilon/2,\delta}} \exp \Big( \sum_{i \leq M} \sum_{j \leq n} \vo_i(j) z_i^j(\bvs) \Big) \, d \nu_{M}^n(\bvt) \Big \rangle_{M,N}.
	\end{align*}
	\textit{Step 4:} Combining the inequalities in Step 1, Step 2, and Step 3 with the factorization \eqref{eq:ass_split} finishes the proof. 
\end{proof}

We now derive a lower bound for the second term appearing in \eqref{eq:ass_step_2}.

\begin{lem}\label{lem:ASSsecondterm}
	We have
	\begin{equation*}
	\lim_{N \to \infty} -\E \log \frac{Z_{N}(\bQ,\epsilon)}{Z_{M,N}(\bQ,\epsilon)} \geq \lim_{N \to \infty} - \E \log \Big\langle  \exp \sqrt{M} y(\bvs) \Big\rangle_{M,N}.
	\end{equation*}
\end{lem}

\begin{proof}
	The proof by Gaussian interpolation is standard, for example \cite[Lemma 2]{CASS}. Consider the interpolating Hamiltonian,
	\begin{equation*}
	H_t(\bvs) = \sum_{j \leq n}  \Big( H_{M,N}^j(\bvs) + \sqrt{t} \sqrt{M} Y^j(\bvs)  + \sqrt{1 - t} \sqrt{M}  y^j(\bvs) \big) \Big),
	\end{equation*}
	Recalling \eqref{eq:cavity}, consider the corresponding interpolating function,
	\begin{align*}
	\phi(t) &= \E \log \int_{Q_N^\epsilon} \exp H_t(\bvs) d \lambda_N(\bvs).
	\end{align*} 
	Differentiating $\phi$, we have
	\begin{equation*}
	\phi'(t) =  \frac{1}{2} \E \bigg\langle \E \frac{\partial H_{t} (\bvs^1)}{ \partial t} \cdot H_{t}(\bvs^1) - \E \frac{\partial H_{t} (\bvp^{\,1})}{ \partial t} \cdot H_{t}(\bvs^2) \bigg\rangle_t 
	\end{equation*} 
	where $\langle \cdot \rangle_t$ is the Gibbs average on $Q_N^\epsilon$ with respect to the Hamiltonian $H_t(\bvs)$. The covariances are given by
	\begin{equation*}
	\E \frac{\partial H_{t} (\bvs^1)}{ \partial t} \cdot H_{t}(\bvs^2)  = M\sum_{j,j' \leq n} \Big( \E Y^j(\bvs^1) Y^{j'}(\bvs^2) - \E y^j(\bvs^1) y^{j'}(\bvs^2) \Big) = \sO\Big(\frac{M}{N}\Big),
	\end{equation*}
	for any $\bvs^1, \bvs^2 \in Q_N^\epsilon$. Integrating $\phi'(t)$, we have 
	\begin{align*}
	\phi(1) = \phi(0) +  \sO\Big(\frac{M}{N}\Big).
	\end{align*}
	Notice \eqref{eq:cavity} implies $\phi(1) = Z_N(\bQ,\epsilon)$. Taking $N \to \infty$, and normalizing both sides by $Z_{M,N}(\bQ,\epsilon)$ finishes the proof.
\end{proof}

The proof of Lemma \ref{lem:ass_unpert} follows by applying Lemma \ref{lem:ASSfirstterm} and Lemma \ref{lem:ASSsecondterm} to \eqref{eq:ass_step_2}. In particular, we have the lower limit of \eqref{eq:ass_step_2} is bounded below by
	\begin{align*}
	&\frac{1}{M} \liminf_{N \to \infty} \bigg( \E \log \Big \langle \int_{\Omega_M^{\epsilon / 2,\delta}} \exp \Big( \sum_{i \leq M} \sum_{j \leq n} \vo_i(j) z_i^j(\bvs) \Big) \, d \nu_{M}^n(\bvo) \Big \rangle_{M,N} - \E \log \Big\langle  \exp \sqrt{M} y(\bvs) \Big\rangle_{M,N} \Big) \bigg) \notag\\
	&\qquad -  L \delta,
	\end{align*}
finishing the proof of Lemma \ref{lem:ass_unpert}.

\section{Perturbation, Ghirlanda--Guerra identities, and their consequences}\label{sec:modpert}

Using the Aizenman--Sims--Starr scheme, we can approximate the lower bound of the free energy with continuous functionals of the distribution of the overlap array. In particular, we have the terms 
\begin{equation*}
\E \log \Big \langle \int_{\Omega_M^{\epsilon / 2, \delta}} \exp \Big( \sum_{i \leq M} \sum_{j \leq n} \vo_i(j) z_i^j(\bvs) \Big) \, d \nu_{M}^n(\bvo) \Big \rangle_{M,N} \text{ and } \E \log \Big\langle  \exp \sqrt{M} y(\bvs) \Big\rangle_{M,N} ,
\end{equation*}
appearing in Lemma \ref{lem:ass_unpert} are continuous functionals of the distributions of the overlap array $(R_{\ell,\ell'})_{\ell,\ell' \geq 1}$ under the Gibbs measure $\E (G_N)^{\otimes \infty}$ \cite[Theorem 1.3]{PPotts}. Before computing the value of the lower bound in the limit, we must first understand the limiting distribution of this overlap. Our main tool is a perturbation of the Gibbs measure that, in the limit, will force the overlaps to satisfy the matrix version of the Ghirlanda--Guerra identities \cite[Theorem 3]{PVS} that in turn imply a powerful synchronization property \cite[Theorem 4]{PVS} in addition to the main consequences of the usual identities \cite[Section 3]{PBook}. These consequences will be summarized at the end of this section.

In this section, we introduce this perturbation of the Hamiltonian. We face two main obstacles. Firstly, the usual proof of the Ghirlanda--Guerra identities requires the self-overlaps $\bR(\bvs,\bvs)$ to be constant, which is not immediate in our setting because self-overlaps are only constrained to lie within an $\epsilon$ window $\bQ$. Secondly, we need to find a suitable perturbation to give us the matrix version of the Ghirlanda--Guerra identities. Both of these issues are resolved in detail in Section 4 and Section 5 of \cite{PVS}. They can be adapted to our setting with a few minor modifications. 
\subsection{Modified Coordinates:} We begin by introducing a transformation of the coordinates that was used to control the self overlaps in the vector spin models \cite[Section 3]{PVS}. This transformation will fix the self overlaps allowing us to apply the usual proof of the Ghirlanda--Guerra identities.  

We use essentially the same change of variables as defined in \cite[Section 3]{PVS} with two main differences. Firstly, since we only need to find a bound for positive definite constraints $\bQ$, we do not need to truncate the constraints like in \cite{PVS}. Secondly, the spins $\vs_i$ are bounded by a universal constant in the vector spin models, while the individual spins in the spherical models have entries bounded by $N$. In our setting, we will need to use a slightly different approach to obtain the relevant bounds on the distortion. 

Let $\lambda_{min} (\bQ) > 0$ denote the smallest eigenvalue of $\bQ$. We first state this transformation as it appears in Section 3 of \cite{PVS}. 
\begin{lem}\label{lem:modifiedcoords}
	\cite[Lemma 4]{PVS} Let $\epsilon < \lambda_{min} (\bQ)$. For each positive definite matrix $\bR$ such that $\| \bR - \bQ \|_\infty \leq \epsilon$, there exists a positive semidefinite matrix $\bA = \bA(\bR)$ such that $\bA \bR \bA^\trans = \bQ$.
	
	Furthermore, we have the bounds
	\begin{equation}\label{eq:bound_A}
	\tr \Big( (\bA - \bI)\bR (\bA - \bI)^\trans \Big) \leq L \sqrt{\epsilon}
	\end{equation}
	and, for any $\bR_1$, $\bR_2$ such that both $\| \bR_1 - \bQ \|_\infty \leq \epsilon$ and $\| \bR_2 - \bQ \|_\infty \leq \epsilon$,
	\begin{equation}\label{eq:bound_B}
	\| \bA(\bR_1) - \bA(\bR_2) \|_\infty \leq \frac{L}{\epsilon} \| \bR_1 - \bR_2 \|_\infty.
	\end{equation}
\end{lem}
In the spherical model, we will also need uniform control on $\|\bA(\bR)\|_\infty$. Since our constant $\bQ$ is positive definite, this fact follows as an immediate consequence of \eqref{eq:bound_A} and \eqref{eq:bound_B}.
\begin{cor}\label{cor:boundc}
	If $\epsilon \leq 1$, each matrix $A(\bR)$ constructed in Lemma \ref{lem:modifiedcoords}, also satisfies the bound
	\begin{equation}\label{eq:bound_C}
	\|A(\bR)\|_\infty \leq L.
	\end{equation}
\end{cor}

\begin{proof}
	We first find a bound on $\bA = \bA(\bQ)$. Since $\bQ$ is positive definite, by the Cholesky decomposition, there exists an invertible matrix $\bB$ such that $\bQ = \bB\bB^\trans$. By \eqref{eq:bound_A}, we have
	\begin{equation*}
	\tr \Bigl( (\bA - \bI)\bQ (\bA - \bI)^\trans \Bigr) = \tr \Bigl( (\bA - \bI) \bB \bB^\trans (\bA - \bI)^\trans \Bigr) = \| (\bA - \bI) \bB \|_{F}^2 \leq L \sqrt{\epsilon}.
	\end{equation*}
	By norm equivalence, we see
	\begin{equation*}
	\| (\bA - \bI) \bB \|_2 \leq \| (\bA - \bI) \bB \|_{F}^2 \leq L \sqrt{\epsilon}.
	\end{equation*}
	Since $\bB$ is invertible and the $\| \cdot \|_2$ norm is sub-multiplicative we have,
	\begin{equation*}
	\| \bA - \bI \|_2 = \| (\bA - \bI) \bB \bB^{-1} \|_2 \leq \| (\bA - \bI) \bB \|_2 \|\bB^{-1}\|_2.
	\end{equation*}
	Therefore, by norm equivalence, if we assume $\epsilon \leq 1$,
	\begin{equation*}
	\| \bA - \bI \|_\infty \leq \sqrt{n} \| \bA - \bI \|_2  \leq L \sqrt{n} \sqrt{\epsilon} \|\bB^{-1}\|_2 \leq L \sqrt{n} \|\bB^{-1}\|_2,
	\end{equation*}
	which implies $\bA(\bQ)$ is uniformly bounded for all $\epsilon \leq 1$.
	
	Furthermore, by \eqref{eq:bound_B}, for any $\bA(\bR)$ such that $\|\bR - \bQ\|_\infty \leq \epsilon$, we have
	\begin{equation}\label{eq:bound_D}
	\| \bA(\bR) - \bA(\bQ) \|_\infty \leq \frac{L}{\epsilon} \| \bR - \bQ \|_\infty \leq L.
	\end{equation}
	Therefore, all matrices $\bA(\bR)$ lie within a closed ball around $\bA(\bQ)$, which implies that $\| \bA(\bR) \|_\infty$ is uniformly bounded for all $\epsilon \leq 1$.
\end{proof}

\noindent\textbf{Remark:} Corollary~\ref{cor:boundc} also holds if we assume is $\bQ$ is only positive semidefinite. The matrix $\bA(\bR)$ has an explicit construction in the proof of Lemma 4 in \cite{PVS}, that only depended on a subset of the eigenvalues of $\bQ$. Therefore, there are finitely many possible constructions of $\bA(\bQ)$, so we can apply the bound \eqref{eq:bound_D} to each possible values $\bA(\bQ)$ to conclude the uniform bound $\|A(\bR)\|_\infty \leq L$.
\\

Lemma~\ref{lem:modifiedcoords} implies there exists a coordinate transform that fixes the self overlaps.  For each $\bvs \in Q_N^\epsilon$, suppose $\bA_\sigma = \bA(\bR(\bvs,\bvs))$ is chosen as in Lemma~\ref{lem:modifiedcoords}. Denote the modified coordinates by $\bhs = ( A_{\sigma} \sigma_i )_{i \leq n} := \bA_\sigma \bvs$ and observe the corresponding modified overlap satisfies
\begin{equation}\label{eq:modifiedoverlaps}
\bR(\bhs,\bhs) = \bR(\bA_{\sigma} \bvs,\bA_{\sigma} \bvs) = \frac{1}{N} \sum_{i \leq N} (\bA_{\sigma} \sigma_i) (\bA_{\sigma} \sigma_i)^\trans = \bA \bR(\bvs,\bvs)  \bA^\trans = \bQ.
\end{equation}

The bounds \eqref{eq:bound_A}, \eqref{eq:bound_B}, and \eqref{eq:bound_C} are used to show the modified overlap matrix is close to the usual overlap. Notice that,
\begin{align*}
\| \bR( \bhs^\ell, \bhs^{\ell'} ) - \bR( \bvs^\ell, \bvs^{\ell'} ) \|_\infty \leq  \| \bR( \bhs^\ell, \bhs^{\ell'} ) - \bR( \bvs^\ell, \bhs^{\ell'} ) \|_\infty  +  \| \bR( \bvs^\ell, \bhs^{\ell'} ) - \bR( \bvs^\ell, \bvs^{\ell'} ) \|_\infty.
\end{align*}
To control the first term, by the Cauchy--Schwarz inequality we have,
\begin{align*}
\| \bR( \bhs^\ell, \bhs^{\ell'} ) - \bR( \bvs^\ell, \bhs^{\ell'} ) \|_\infty &\leq \sup_{j,j' \leq n} \frac{1}{N} \bigg| \sum_{i = 1}^N \bA_{\sigma}\vs^\ell_i(j) \bA_{\sigma} \bvs^{\ell'}_i(j') - \sigma^{\ell}_i(j) \bA_{\sigma}\bvs^{\ell'}_i(j') \bigg|
\\&\leq  \frac{1}{N} \sup_{j,j' \leq n} \|(\bA_{\sigma} - \bI)\bvs^{\ell}(j)\| \|\bA_{\sigma} \bvs^{\ell'}(j')\| 
\\&\leq \sup_{j \leq n}  \frac{\|(\bA_{\sigma} - \bI)\bvs^{\ell}(j)\|}{\sqrt{N}} \| \bA_{\sigma}\|_\infty
\\&\leq \|\bA_{\sigma} \|_\infty \tr ( \bR( (\bA_{\sigma} - \bI) \bvs^{\ell}, (\bA_{\sigma} - \bI) \bvs^{\ell}) )^{1/2}.
\end{align*}
Using observation \eqref{eq:modifiedoverlaps}, the bounds \eqref{eq:bound_A} and \eqref{eq:bound_C} imply
\begin{equation*}
\| \bR( \bhs^\ell, \bhs^{\ell'} ) - \bR( \bvs^\ell, \bhs^{\ell'} ) \|_\infty \leq L \epsilon^{1/4}.
\end{equation*}
A similar computation applied to the second term gives a similar bound,
\begin{equation*}
\| \bR( \bvs^\ell, \bhs^{\ell'} ) - \bR( \bvs^\ell, \bvs^{\ell'} ) \|_\infty \leq L \epsilon^{1/4}.
\end{equation*}
Therefore, the modified overlap only differs from the overlap by a factor of $\epsilon^{1/4}$,
\begin{equation}\label{eq:overlaps_bound}
\| \bR( \bhs^\ell, \bhs^{\ell'} ) - \bR( \bvs^\ell, \bvs^{\ell'} ) \|_\infty \leq L \epsilon^{1/4}.
\end{equation}
The bounds \eqref{eq:overlaps_bound} and \eqref{eq:bound_B} will ensure this change of variables will not affect the limiting values in the perturbed Aizenman--Sims--Starr scheme that we introduce next.
\subsection{Perturbed Hamiltonian:} We now define the perturbation that will force the overlaps to satisfy the matrix version of the Ghirlanda--Guerra identities in \cite{PVS}. This perturbation is identical to the one introduced in Section 5 of \cite{PVS}. We summarize the key steps below.

We denote the family of parameters
\begin{equation}\label{SHOWLABEL}
\theta = (p,m,n_1, \dots, n_m, \vnu^{\, 1}, \dots, \vnu^{\, m}).
\end{equation}
For each $\theta$, there exists Gaussian processes $h_\theta(\bvs)$ indexed by $\bvs \in  S_N^n$ with mean $0$ and covariance
\begin{equation}\label{SHOWLABEL}
C_{\ell, \ell'}^\theta = \Cov\bigl(h_\theta(\bvs^\ell), h_\theta(\bvs^{\ell'}) \bigr) = \prod_{j \leq m} \big( R^{\circ p}_{\ell, \ell'} \vnu^{\, j}, \vnu^{\,j} \big)^{n_j}.
\end{equation}
Furthermore, for $\nu \in [-1,1]^n$ and $\bvs \in S_N^n$, the covariance is bounded by $n^{2p(n_1 + \dots + n_m)}$. We denote the countable set of parameters with
\begin{equation}\label{SHOWLABEL}
\Theta = \{ \theta \mmm p \geq 1, m \geq 1, n_1, \dots, n_m \geq 1, \vnu^{\, 1}, \dots , \vnu^{\, m} \in ( [-1,1] \cap \mathbb{Q} )^n \}.
\end{equation} 
Let $j_0:  ([-1,1] \cap \mathbb{Q})^n \to \N$ be a one-to-one function. We denote an enumeration of $\theta \in \Theta$ with
\begin{equation}\label{SHOWLABEL}
j(\theta) = p + n_1 + \dots + n_m + j_0(\vnu^{\, 1}) + \dots + j_0(\vnu^{\, m}) + 22m.
\end{equation}
Let $(u_\theta)_{\theta \in \Theta}$ be a random sequence of \iid uniform random variables in $[1,2]$. We define the interpolating Hamiltonian,
\begin{equation}\label{SHOWLABEL}
h_N(\bvs) = \sum_{\theta \in \Theta} 2^{-j(\theta)} n^{2(n_1 + \dots+ n_m)}  u_\theta h_\theta(\bvs).
\end{equation}
The covariance of this process is bounded by $1$, and given explicitly by
\begin{equation}\label{eq:covfunc}
\Cov\bigl( h_N(\bvs^\ell), h_N(\bvs^{\ell'}) \bigr) = \sum_{\theta \in \Theta} 2^{-2j(\theta)} n^{4(n_1 + \dots+ n_m)}  u^2_\theta 
\prod_{j \leq m} \big( R^{\circ p}_{\ell, \ell'} \vnu^{\,j}, \vnu^{\,j} \big)^{n_j}.
\end{equation} 
For $\frac{1}{4} < \gamma < \frac{1}{2}$, we denote the sequence $s_N = N^\gamma$. Recall the modified coordinates defined in the previous section denoted with $\bhs = (\bA_{\sigma} \vs_i)_{i \leq N}$. We define the perturbed Hamiltonian
\begin{equation}\label{SHOWLABEL}
H_N^{pert}(\bvs) = H_N(\bvs) + s_N h_N(\bhs),
\end{equation}
and the corresponding perturbed partition function
\begin{equation}\label{SHOWLABEL}
Z_N^{pert}(\bQ, \epsilon) = \int_{Q_N^\epsilon} \exp \Bigl( H_N^{pert}(\bvs) + \sum_{i \leq N} \sum_{j \leq n}  \vh(j) \vs_i(j) \Bigr) \, d \lambda_N^n(\bvs).
\end{equation}
Since $\frac{s_N^2}{N} \to 0$, a straightforward Gaussian interpolation argument shows
\begin{equation}\label{eq:pert_freenergy}
\liminf_{N \to \infty} \frac{1}{N} \E \log Z_N(\bQ,\epsilon) = \liminf_{N \to \infty} \frac{1}{N} \E \log Z_N^{pert}(\bQ,\epsilon).
\end{equation}
\subsection{Perturbed Aizenman--Sims--Starr Scheme:} The Aizenman--Sims--Starr scheme proved in Section~\ref{sec:ASS} has to be modified slightly to account for the extra perturbation term in the Hamiltonian. Let $\langle \cdot \rangle_{pert}$ be the average on $Q_N^\epsilon$ with respect to the Gibbs measure
\begin{equation}\label{eq:Gpert}
G_N^{pert}(\bvs) = \frac{\exp\Big(  H_N^{pert}(\bvs) + \sum_{i \leq N} \vh(j) \vs_i(j) \Big)}{Z_N^{pert}(\bQ, \epsilon)}.
\end{equation}
The following modification of Lemma \ref{lem:ass_unpert} will be used in the proof of lower bound.
\begin{lem} \label{lem:ASS_Pert}
	For $s_N = N^\gamma$, $\vh = \vec{0}$ and $\bhs = (\bA_{\sigma} \vs_i)_{i \leq N}$ we have
	\begin{align}
	\liminf_{N \to \infty} \frac{1}{N} \E \log Z_N^{pert} &\geq \frac{1}{M} \liminf_{N \to \infty} \bigg( \E \log \Big \langle \int_{\Omega_M^{\epsilon/2,\delta}} \exp \Big( \sum_{i \leq M} \sum_{j \leq n} \vo_i(j)  z_i^j(\bhs)  \Big) \, d \nu_{M}^n(\bvo) \Big \rangle_{pert} \notag
	\\&\qquad - \E \log \Big\langle  \exp \sqrt{M} y(\bhs) \Big\rangle_{pert} \bigg) - L \delta  - L \epsilon^{1/4}. \label{eq:ASS_Pert}
	\end{align}
\end{lem}

\begin{proof}
	Only a small modification needs to be made to adapt the proof of Lemma \ref{lem:ass_unpert} to this setting. We can leave the perturbation term $s_N h_N$ out of the interpolation in the proof of Lemma \ref{lem:ASSsecondterm} and keep the rest of the proof unchanged. To adapt Lemma \ref{lem:ASSfirstterm}, we have to control
	\begin{align}\label{eq:ass_split_pert}
	\E \log \dfrac{Z^{pert}_{M + N}(\bQ,\epsilon)}{\int_{\Omega_M^{\epsilon/2,\delta}} J^{pert}_{M,N} \, d{\nu}_{M,N}^n(\bvo) } + \E \log \dfrac{\int_{\Omega_M^{\epsilon/2,\delta}} J^{pert}_{M,N} \, d{\nu}_{M,N}^n(\bvo) }{\int_{\Omega_M^{\epsilon/2,\delta}} J^{pert}_{M,N} \, d{\nu}_{M}^n(\bvo) } + \E \log \dfrac{\int_{\Omega_M^{\epsilon/2,\delta}} J^{pert}_{M,N} \, d{\nu}_{M}^n(\bvo) }{Z^{pert}_{M,N}(\bQ,\epsilon)},
	\end{align}
	where
	\begin{align*}
	J^{pert}_{M,N} = \int_{Q_N^{\epsilon}} \exp \Big(\sum_{j \leq n} H^j_{M,N} (\bvs) + \sum_{j \leq n} \sum_{i \leq M} \vo_i(j) z_i^j(\bvs) + s_N h_N (\bhs) \Big) \, d \lambda^n_N(\bvs).
	\end{align*}
	The perturbation term appears as $s_{N + M} h_{N + M} ( \bA_{\rho}\ubvp)$ in $Z^{pert}_{M + N}(\bQ,\epsilon)$, but we need it to appear as $s_N h_N (\bhs)$ to match the normalization. This issue is resolved by reproving the bound,
	\begin{align*}
	\liminf_{N \to \infty} \E \log \dfrac{Z^{pert}_{M + N}(\bQ,\epsilon)}{\int_{\Omega_M^{\epsilon/2,\delta}} J^{pert}_{M,N} \, d{\nu}_{M,N}^n(\bvo) } \geq - LM\delta.
	\end{align*}
	Consider the interpolating Hamiltonian
	\begin{equation}\label{eq:triplehamiltonian_pert}
	H^{pert}_t(\bvp) = \sum_{j \leq n}H^j_{1,t}(\bvp\,) + H^j_{2,t} (\bvp\,) + H^j_{3,t} (\bvp\,) + H_{4,t} (\bvp\,).
	\end{equation}
	where the Hamiltonians are defined in \eqref{eq:ham_1}, \eqref{eq:ham_2}, \eqref{eq:ham_3} and
	\begin{equation*}
	H_{4,t}(\bvp) = \sqrt{t} s_{N + M} h_{N + M} \big( \bA_{\rho} \ubvp \big) + \sqrt{1 - t} s_N h_N (\bhs),
	\end{equation*}
	
	After applying Gaussian integration by parts, we will need to control
	\begin{align}\label{eq:covpert}
	\bigg| \E \frac{d H_t(\bvp^1)}{dt} H_t(\bvp^2) \bigg| &= \Big| s_{N + M}^2 \E h_{N + M} \big( \bA_{\rho^1} \ubvp^1 ) h_{N + M} \big( \bA_{\rho^2}  \ubvp^2  ) -  s_N^2 \E h_N (\bhs^1) h_N (\bhs^2) \Big| \notag
	\\&= \Big| (N + M)^{2 \gamma} g(\bR(\bA_{\rho^1} \ubvp^1,\bA_{\rho^2} \ubvp^2  )) -  N^{2 \gamma} g(\bR(\bA_{\sigma^1}\bvs^1,\bA_{\sigma^2}\bvs^2)) \Big|,
	\end{align}
	where $g$ is the covariance function of $h_N$ given by \eqref{eq:covfunc}. The function $g$ and its derivatives is bounded on compacts uniformly for all parameters $u_\theta$. Using \eqref{eq:overlaprhosplit} and \eqref{eq:limitcoeffN}
	\begin{align*}
	&\| \bR(\ubvp^{\,1}, \ubvp^{\,2}) - \bR(\bvs^1,\bvs^2) \|_\infty 
	\\&= \sup_{j,j' \leq n} \bigg| \frac{(a^{j}_{M + 1}(\bo^1)a^{j'}_{M + 1}(\bo^2) N - N) R^{j,j'}(\bvs^1, \bvs^2)  }{M + N} - \frac{MR^{j,j'}(\bvs^1, \bvs^2)}{M + N} + \frac{M R^{j,j'}(\bto^1, \bto^2)  }{M + N}  \bigg|
	\\&= \sO(N^{-1}),
	\end{align*}
	and therefore, by Lemma \ref{lem:modifiedcoords},
	\begin{align*}
	\| \bR(\bA_{\rho^1} \ubvp^1,\bA_{\rho^2}\ubvp^2  ) - \bR(\bA_{\sigma^1}(\bvs^1),\bA_{\sigma^2}(\bvs^2)) \|_\infty&=	\| \bA_{\rho^1} \bR(\ubvp^1, \ubvp^2)\bA_{\rho^2}^\trans - \bA_{\sigma^1} \bR(\bvs^1,\bvs^2) \bA_{\sigma^2}^\trans \|_\infty 
	\\&\leq \frac{L	\| \bR(\ubvp^{1}, \ubvp^{2}) - \bR(\bvs^1,\bvs^2) \|_\infty }{\epsilon} 
	\\&= \sO((N\epsilon)^{-1}).
	\end{align*}
	Using the Taylor series of $g(\bR(\bA_{\rho^1}\ubvp^1,\bA_{\rho^2}\ubvp^2)$ around $\bR(\bA_{\sigma^1}\bvs^1,\bA_{\sigma^2}\bvs^2)$, we see 
	\begin{equation*}
	(N + M)^{2 \gamma} g(\bR(\bA_{\rho^1}\ubvp^1,\bA_{\rho^2}\ubvp^2)) = (N + M)^{2 \gamma} g(\bR(\bA_{\sigma^1}\bvs^1,\bA_{\sigma^2}\bvs^2) ) + \sO( N^{-1 - 2\gamma}/\epsilon ).
	\end{equation*}
	Since $(N + M)^{2 \gamma} - N^\gamma = \sO(N^{-1 - 2\gamma})$ we see
	\eqref{eq:covpert} is
	\begin{align*}
	\Big| (N + M)^{2 \gamma} g(\bR(\bA_{\rho^1} \ubvp^1,\bA_{\rho^2} \ubvp^2)) -  N^{2 \gamma} g(\bR(\bA_{\sigma^1}\bvs^1,\bA_{\sigma^2}\bvs^2)) \Big| = \sO( N^{- (1 - 2\gamma)}/\epsilon ).
	\end{align*}
	The above bound holds uniformly for $u_\theta$, so combined with the fact $\gamma < \frac{1}{2}$, this means that the replacing $s_{N + M} h_{N + M}( \bA_{\rho} \ubvp)$ with $s_N h_N (\bhs)$ is small enough to not change the errors in the proof of Lemma \ref{lem:ass_unpert}. The remainder of the proof of Lemma \ref{lem:ASSfirstterm} is unchanged, we arrive at
	\begin{align}
	\liminf_{N \to \infty} \frac{1}{N} \E \log Z_N^{pert} &\geq \frac{1}{M} \liminf_{N \to \infty} \bigg( \E \log \Big \langle \int_{\Omega_M^{\epsilon/2,\delta}} \exp \Big( \sum_{i \leq M} \sum_{j \leq n} \vo_i(j)  z_i^j(\bvs)  \Big) \, d \nu_{M}^n(\bvo) \Big \rangle_{pert} \notag
	\\&\qquad - \E \log \Big\langle  \exp \sqrt{M} y(\bvs) \Big\rangle_{pert} \bigg) - L \delta \label{eq:ass_pert_nomod}.
	\end{align}
	When we characterize the limiting distribution of the overlap array, we will require the self overlaps to be constant. Replacing $\bvs$ with the modified coordinates $\bhs$ in the cavity fields achieves this. Starting from \eqref{eq:ass_pert_nomod}, an interpolation argument will prove that the cavity fields can be replaced with
	\begin{align}\label{eq:pert_z}
	\E \log \Big \langle \int_{\Omega_M^{\epsilon/2,\delta}} \exp \Big( \sum_{i \leq M} \sum_{j \leq n} \vo_i(j) z_i^j(\bhs) \big) \Big) \, d \nu_{M}^n(\bvo) \Big \rangle_{pert},
	\end{align}
	and
	\begin{align}\label{eq:pert_y}
	\E \log \Big\langle  \exp \sqrt{M} y(\bhs) \Big\rangle_{pert}
	\end{align}
	at the cost of $L \epsilon^{1/4}$ error. We only prove \eqref{eq:pert_z} because the proof of \eqref{eq:pert_y} is almost identical. 
	
	Consider the Hamiltonian,
	\begin{equation*}
	Z_i^j(\bvs;t) = \sqrt{t} z_i^j(\bvs) + \sqrt{1 - t} z_i^j(\bhs),
	\end{equation*}
	and the corresponding interpolating function
	\begin{equation*}
	\phi(t) = \E \log \Big \langle \int_{Q_M^{\epsilon/2}} \exp \Big( \sum_{i \leq M} \sum_{j \leq n} \vt_i(j) Z_i^j(\bvs;t) \Big) \, d \lambda_{M}^n(\bvt) \Big \rangle_{pert}.
	\end{equation*}
	Let $\hat{\bR}_{\ell,\ell'} := R(\bhs^\ell, \bhs^{\ell'})$, a standard integration by parts computation will show
	\begin{align*}
	|\phi'(t)| &\leq \Big\| \bR_{1,1} \odot  \Big( \bxi'(\bR_{1,1}) - \bxi'(\hat{\bR}_{1,1}) \Big) - \bR_{1,2} \odot \Big( \bxi'(\bR_{1,2}) - \bxi'(\hat{\bR}_{1,2}) \Big) \Big\|_\infty\\
	&\leq  n^2 \xi'(1) \Big( \|\bR_{1,1} - \hat{\bR}_{1,1}\|_\infty + \|\bR_{1,2} - \hat{\bR}_{1,2}\|_\infty \Big)\\
	&\leq  L \epsilon^{1/4}
	\end{align*}
	since $\|\bR_{1,1} - \hat{\bR}_{1,1}\|_\infty \leq L\epsilon^{1/4}$ by Lemma \ref{lem:modifiedcoords} and \eqref{eq:overlaps_bound}. Integrating the quantity above implies
	\begin{equation*}
	\phi(0) \geq \phi(1) - \sup_{t \in [0,1]} |\phi'(t)| \geq \phi(1) - L\epsilon^{1/4}.
	\end{equation*}
	The bound for \eqref{eq:pert_y} is similar to above, and is proved using the interpolation
	\begin{equation*}
	Y(\bvs;t) = \sqrt{t} y(\bvs) + \sqrt{1 - t} y(\bhs).
	\end{equation*}
	Applying the bounds \eqref{eq:pert_z} and \eqref{eq:pert_y} to \eqref{eq:ass_pert_nomod} finishes the proof.
\end{proof}

\noindent\textbf{Remark:} We assumed $\veh = \vec{0}$ in the computations above to simplify notation. If $\veh$ was non-zero, then the lower bound \eqref{eq:ASS_Pert} in Lemma \ref{lem:ASS_Pert} is of the form
\begin{align}
&\frac{1}{M} \liminf_{N \to \infty} \bigg( \E \log \Big \langle \int_{\Omega_M^{\epsilon/2,\delta}} \exp \Big( \sum_{i \leq M} \sum_{j \leq n} \vo_i(j) \big( z_i^j(\bhs) + \vh(j) \big) \Big) \, d \nu_{M}^n(\bvo) \Big \rangle_{pert} \label{eq:z_externalfield}
\\&\qquad - \E \log \Big\langle  \exp \sqrt{M} y(\bhs) \Big\rangle_{pert} \bigg) - L\delta \label{eq:h_externalfield}.
\end{align}
where $\langle \cdot \rangle_{pert}$ is the average on $Q_N^\epsilon$ with respect to the Gibbs measure with external field,
\begin{equation}\label{eq:gibbs_field}
G_N(d\bvs) \propto \exp\Big(  H_{M,N}(\bhs) + \sum_{i \leq N} \sum_{j \leq n} \vh(j) \vs_i(j) \Big) \, d\lambda_N^n(\bvs).
\end{equation}

The bound \eqref{eq:h_externalfield} follows by a simple modification of the above proof. The external field can be decoupled into its cavity and non-cavity coordinates immediately,
\begin{equation}\label{eq:assfield1}
\sum_{i \leq M + N} \sum_{j \leq n}  \vh(j)\vp_i(j) = \sum_{i \leq N} \sum_{j \leq n}  \vh(j)\vs_i(j) + \sum_{i \leq M} \sum_{j \leq n}  \vh(j)\vo_i(j).
\end{equation}
The first summation appears in the Gibbs average \eqref{eq:gibbs_field} and the second summation appears in the cavity field term \eqref{eq:z_externalfield}. However, the external field in the exponent of \eqref{eq:step1_sub} will appear as
\begin{equation}\label{eq:assfield2}
\sum_{i \leq N} \sum_{j \leq n}  \vh(j) a_{M + 1}^j \vs_i(j) + \sum_{i \leq M} \sum_{j \leq n}  \vh(j) a_{i}^j \vo_i(j)
\end{equation}
in Step 1 of the proof of Lemma \ref{lem:ASSfirstterm}. To resolve this issue, notice for $\bvo \in \Omega_M^{\epsilon/2,\delta}$ each term $a_{\ell}^j(\bvo) \to 1$ uniformly on $\Omega_M^{\epsilon/2,\delta}$ for all $\ell \leq M$. For the $M + 1$ coefficient, we also have
\begin{equation*}
\lim_{N \to \infty } N|a^j_{M + 1} - 1| =  \bigg( \frac{M}{2} - \frac{\| \bvo(j) \|^2}{2} \bigg) \leq \frac{M \delta}{2} .
\end{equation*}
Therefore, by the Cauchy--Schwarz inequality, for all $(\bvs,\bvo) \in Q_N^\epsilon \times \Omega_M^{\epsilon/2,\delta}$ and $j \leq n$ we have
\begin{align*}
& \bigg\| \sum_{i \leq N} \vh(j) a_{M + 1}^j \vs_i(j) + \sum_{i \leq M}  \vh(j) a_{i}^j \vo_i(j) - \sum_{i \leq N} \vh(j) \vs_i(j) - \sum_{i \leq M} \vh(j) \vo_i(j) \bigg\|_\infty
\\&\leq \| \veh \|_\infty N |a^j_{M + 1} - 1| + M \|\veh\|_\infty \sup_{i \leq M} | a_i^j - 1|
\\&\leq L M \delta 
\end{align*}
for $N$ sufficiently large. Therefore, we can replace the external field in \eqref{eq:assfield2} with \eqref{eq:assfield1} and absorb the $LM\delta$ error into the right hand side of \eqref{eq:ass_pert_nomod}. 
\subsection{Consequences of the Perturbation:} The lower bound \eqref{eq:h_externalfield} is a continuous functional of the distribution of the modified arrays $(\bR(\bhs^\ell,\bhs^{\ell'}))_{\ell,\ell' \geq 1}$ under the Gibbs average $\E(G_N^{pert})^{\otimes \infty}$ \cite[Lemma 8]{PPotts}, so it suffices to study the distribution of the modified array. To this end, we state matrix version of the Ghirlanda--Guerra identities and several of its consequences. These are identical to \cite[Section 5]{PVS} and can now be applied in this setting with no modification.

The entries of the overlaps are in $[-1,1]$, so the probability distributions on finite dimensional subsets of the infinite array are tight. Therefore, by the selection theorem, there exists a subsequence such that all finite dimensional distributions of $(\bR(\bhs^\ell,\bhs^{\ell'}))_{\ell,\ell' \geq 1}$ converge weakly. Furthermore, there exists a non-random sequence of parameters $(u^N_\theta)$ (see \cite[Lemma 5]{PVS} and \cite[Lemma 3.3]{PBook}), possibly changing in $N$, such that the limiting array, denoted by $(\hat\bR_{\ell,\ell'})_{\ell,\ell' \geq 1}$ also satisfies a matrix version of the Ghirlanda--Guerra identities.

Consider $k$ replica of this limiting array, $\hat R^k = (\hat{\bR}_{\ell, \ell'})_{\ell, \ell' \leq k}$, we have:

\begin{lem}\cite[Theorem 3]{PVS} \label{lem:GGI}
	Given any measurable function $\phi: \R^m \to \R$ and $f = f(R^k)$, the array satisfies the Ghirlanda--Guerra identities
	\begin{equation}\label{eq:GGI}
	\E f(\hat R^k) {\bC}_{1, k + 1} = \frac{1}{k} \E f(\hat R^k) \E {\bC}_{1,2} + \frac{1}{k} \sum_{\ell = 2}^{k} \E f(\hat R^k) {\bC}_{1,\ell},
	\end{equation}
	where
	\begin{equation}\label{SHOWLABEL}
	{\bC}_{ \ell, \ell' } = \phi\Big( (\hat{\bR}^{\circ p}_{\ell, \ell'} \nu^{\,1}, \nu^{\,1}), \dots, (\hat{\bR}^{\circ p}_{\ell, \ell'} \nu^{\,m}, \nu^{\,m})  \Big).
	\end{equation}
\end{lem}

We have two main consequences of Lemma \ref{lem:GGI}. If we take $\vec{\nu}_i = e_i$ the standard basis vectors in $\R^n$, \eqref{eq:GGI} implies the traces of the overlap array, denoted by $(T_{\ell,\ell'})_{\ell,\ell' \geq 1} = ( \tr(\hat{\bR}_{\ell, \ell'}) )_{\ell, \ell' \geq 1}$, satisfy the usual Ghirlanda--Guerra identities,
\begin{equation}\label{eq:TRGGI}
\E f(T^k) g(T_{1,k+1}) = \frac{1}{k} \E f(T^k) \E g(T_{1,2}) + \frac{1}{k} \sum_{\ell = 2}^{k} \E f(T^k) g(T_{1,\ell}), \quad g: \R \to \R.
\end{equation}
where $T^k = ( T_{\ell,\ell'})_{\ell, \ell' \leq k}$ is a sample of $k$ replicas from the array of traces and $g$ is a measurable function. In particular, we are able to apply all the consequences of the standard Ghirlanda--Guerra identities to $(T_{\ell,\ell'})_{\ell,\ell' \geq 1}$.

Furthermore, \eqref{eq:GGI} implies a synchronization property for overlap matrices \cite{PVS,PPotts}:
\begin{lem}\label{lem:synch}
	There is a function $\Phi: \R^+ \to \Gamma_n$ such that
	\begin{equation}\label{eq:synch}
	\btR_{\ell, \ell'} = \Phi\bigl( \tr (\btR_{\ell, \ell'} ) \bigr) \text{ a.s.}
	\end{equation}
	Furthermore, this function is non-decreasing, $\Phi(x_1) \leq \Phi(x_2)$ for all $x_1 \leq x_2$, and Lipshitz continuous, $\| \Phi(x_2) - \Phi(x_1) \|_1 \leq L |x_2 - x_1|$. 
\end{lem} 

Lemma \ref{lem:GGI} and Lemma \ref{lem:synch} will allow us to characterize the distribution of the limiting array in the final step of the proof of the lower bound. 

\section{Lower Bound --- Cavity Computations}\label{sec:lwbd}

We now have the tools to prove the lower limit of the free energy. The remainder of the proof is standard and almost identical to other spin glass models (see Chapter 3 of \cite{PBook} or the proof of the lower bound in \cite{PPotts} and \cite{PVS}). We will summarize the steps and reiterate the importance of the synchronization mechanism.

Let $\bQ$ be a positive definite constraint. Starting from the Aizenman--Sims--Starr scheme \eqref{eq:h_externalfield}, we have $\liminf_{N \to \infty} F^\epsilon_N(\bvb, \bQ)$ is bounded below by
\begin{align}
&\frac{1}{M} \liminf_{N \to \infty} \bigg( \E \log \Big \langle \int_{\Omega_M^{\epsilon/2,\delta}} \exp \Big( \sum_{i \leq M} \sum_{j \leq n} \vo_i(j) \big( z_i^j(\bhs) + \vh(j) \big) \Big) \, d \nu_{M}^n(\bvo) \Big \rangle_{pert} \label{eq:zandyfunc}\\
&\qquad - \E \log \Big\langle  \exp \sqrt{M} y(\bhs) \Big\rangle_{pert} \Big) \bigg) -  L \delta  - L\epsilon^{1/4}. \notag
\end{align}
From \cite[Lemma 8]{PPotts}, the averages on \eqref{eq:zandyfunc} are continuous functionals of the distribution of the modified infinite array $(\hat{\bR}^{N,M}_{\ell,\ell'})_{\ell, \ell' \geq 1} := (\bR(\bhs^\ell, \bhs^{\ell'}))_{\ell, \ell' \geq 1}$. To compute this lower bound explicitly, it suffices to understand the limiting distribution of the array $(\hat{\bR}^{N,M}_{\ell,\ell'})_{\ell, \ell' \geq 1}$ under the perturbed Gibbs measure $\E( G^{pert}_N)^{\otimes \infty}$ defined in \eqref{eq:Gpert} for a deterministic choice of parameters $(u_\theta^N)$ such that Lemma~\ref{lem:GGI} holds. By the selection theorem, there exists a subsequence such that
\begin{equation*}
\big(\btR^{N,M}_{\ell,\ell'}\big)_{\ell, \ell' \geq 1} \stackrel{d}{\to} (\btR_{\ell,\ell'}^M)_{\ell,\ell' \geq 1}.
\end{equation*}
The diagonal elements of this array are constant, so by Lemma \ref{lem:GGI}, the limiting array $(\btR_{\ell,\ell'}^M)_{\ell,\ell'}$ satisfies the generalized Ghirlanda--Guerra identities \eqref{eq:GGI} and the synchronization property \eqref{eq:synch}. In particular, there exists a function $\Phi : [0,1] \to \Gamma_n$ such that
\begin{equation*}
\btR_{\ell,\ell'}^M = \Phi ( \tr(\btR_{\ell,\ell'}^M) )
\end{equation*}
almost surely. Recall that $\Phi$ is non-decreasing and Lipschitz. This allows us to approximate its distribution with a random measure generated by the Ruelle probability cascades.

We begin by characterizing the array $(\tr(\btR_{\ell,\ell'}^M))_{\ell,\ell'\geq 1}$ consisting of the traces of the limiting array. As a consequence of the generalized Ghirlanda--Guerra identities \eqref{eq:TRGGI}, the array of traces also satisfies the usual Ghirlanda--Guerra identities. We denote the distribution of $\tr(\btR_{1,2}^M)$ with
\begin{equation}\label{eq:distlimitingarray}
\mu(q) = \pP \bigl( \tr(\btR_{1,2}^M) \leq q \bigr).
\end{equation}
Following the usual proof of the lower bound (see Chapter 3 of  \cite{PBook}) there exists a sequence of probability distributions $(\mu_k)_{k \geq 1}$ such that $\mu_k \to \mu$ in $L^1$,
\begin{equation*}
\lim_{k \to \infty} \int_0^n | \mu_k(q) - \mu(q) | \, dq = 0. 
\end{equation*}
For each $k$, we can encode the discrete probability measures with a sequences of parameters
\begin{equation}\label{eq:paramtrace}
\begin{linsys}{6}
x_{-1} = 0 & < & x_0 & < &  x_1 & < & \dots & < & x_r & = &  1\\
0  &=&  q_0 & < & q_1 & < & \dots & < & q_r & = & n &=& \tr(\bQ)
\end{linsys}
\end{equation}
such that 
\begin{equation}\label{eq:discretedist}
\mu_k (q) = x_p \text{ for } q_p \leq q < q_{p + 1}.
\end{equation}
Let $(v_\alpha)_{\alpha \in \N^r}$ be the Ruelle probability cascades corresponding to \eqref{eq:paramtrace}. Let $(\alpha^\ell)_{\ell \geq 1}$ be an \iid sample from $\N^r$ according to the weights $(v_\alpha)_{\alpha \in \N^r}$, it follows that the array
\begin{equation*}
(T^k_{\ell,\ell'})_{\ell,\ell' \geq 1} = ( q_{\alpha^\ell \wedge \alpha^{\ell'}} )_{\ell, \ell' \geq 1}
\end{equation*}
also converges to $(\tr(\btR_{\ell,\ell'}^M))_{\ell,\ell'\geq 1}$ by Theorem~2.13 and Theorem~2.17 in \cite{PBook}.

From here, we use the synchronization mechanism to recover a sequence of monotone paths in $\Pi$ that describes the distribution of the limiting overlap matrix array $(\btR^{N,M}_{\ell,\ell'})_{\ell, \ell' \geq 1}$. We define
\begin{equation*}
\bQ^k_{\ell,\ell'} = \Phi(T^k_{\ell,\ell'}),
\end{equation*}
and observe $(\bQ^k_{\ell,\ell'})_{\ell,\ell' \geq 1}$ converges to the distribution of $(\btR_{\ell,\ell'}^M)_{\ell,\ell'\geq 1}$ because $\Phi$ is Lipschitz. It also follows that the discrete path
\begin{equation}\label{SHOWLABEL}
\pi_k(x) = \bQ_k \text{ for } x_{k - 1} < x \leq x_{k} \text{ for $0 \leq k \leq r$},\quad \pi(0) = \bm{0}, \quad \pi(1) = \bQ.
\end{equation}
induced by
\begin{equation}\label{eq:paramQ}
\begin{linsys}{6}
x_{-1} = 0 & < & x_0 & < &  x_1 & < & \dots & < & x_r & = &  1\\
0  &=&  \bQ_0 & < & \bQ_1 & < & \dots & < & \bQ_r & = & \bQ.
\end{linsys}
\end{equation}
where $\bQ_\ell = \Phi(q_\ell)$ for $0 \leq \ell \leq r$ is a discretization of the path associated with the limiting array. To see this, recall \eqref{eq:distlimitingarray} and define
\begin{equation*}
\pi(x) := \Phi(\mu^{-1}(x)) \in \Pi,
\end{equation*}
where $\mu^{-1}: [0,1] \to \R^+$ is the quantile distribution of $\mu$. Similarly, for discrete $\mu_k$ given by \eqref{eq:discretedist}, the paths
\begin{equation*}
\pi_k(x) := \Phi(\mu_k^{-1}(x)) \in \Pi,
\end{equation*}
are a discrete approximation of $\pi$ \cite[Equation (71)]{PPotts},
\begin{equation*}
\mathrm{d}(\pi,\pi_k) = \int_0^1 \| \pi(x) - \pi_k (x) \|_1 \, dx \leq n \int_0^1 | \tr(\pi(x)) - \tr(\pi_k (x))| \, dx = n \int_0^1 | \mu(x) - \mu_{k}(x)| \, dx.
\end{equation*} 
In particular, we have $\mathrm{d}(\pi,\pi_k) \to 0$ as $\mu_k \to \mu$ in $L^1$.

Recall the Gaussian processes $Z_i^j(\alpha)$ and $Y(\alpha)$ defined in the \eqref{eq:cov_Z} and \eqref{eq:cov_Y} and consider the following functionals of the discrete paths associated with the approximating arrays $(\bQ_{\alpha^{\ell} \wedge \alpha^{\ell'}})_{\ell, \ell' \geq 1}$:
\begin{align}
f_M^1(\pi) &= \frac{1}{M} \E \log  \sum_{\alpha \in \N^r} v_\alpha \int_{\Omega_M^{\epsilon / 2,\delta}} \exp \Big( \sum_{i \leq M} \sum_{j \leq n} \vo_i(j) \big( Z_i^j(\alpha) + \vh(j)\big) \Big) \, d \nu_{M}^n(\bvo),\label{eq:func1} \\
f_M^2(\pi) &= \frac{1}{M} \E \log  \sum_{\alpha \in \N^r} v_\alpha  \exp \sqrt{M} Y(\alpha).\label{eq:func2}
\end{align}
The covariances of $Z^j_i(\alpha)$ and $z^j_i(\bhs)$, and $Y(\alpha)$ and $y(\bhs)$ are given by the same functions of arrays so the difference of the functionals \eqref{eq:func1}, \eqref{eq:func2} and the functional appearing in \eqref{eq:zandyfunc} can be approximated by the same continuous bounded function of the array \cite[Lemma 8]{PPotts}. In summary, by choosing a discretization $\mu_k$ close enough to $\mu$ in $L^1$, we can find a corresponding discrete path $\hat\pi_M:= \pi_k$ encoded by the sequences \eqref{eq:paramQ} such that
\begin{equation}\label{eq:lowbd_RPC}
\liminf_{N \to \infty} F^\epsilon_N(\bvb, \bQ) \geq f_M^1(\hat\pi_M) - f_M^2(\hat\pi_M) - L\delta  - L\epsilon^{1/4}.
\end{equation}

The lower bound holds for all $M$, so we can take a sub-sequential limit as $M \to \infty$. However, we cannot apply Lemma \ref{lem:sharpupbd} to compute the lower bound, because the paths $\hat\pi_M$ may change in $M$. To resolve this, notice that by monotonicity of the paths, $\hat{\pi}_M \to \hat{\pi}$ along some subsequence \cite[Section 7]{PPotts}. Furthermore, there exists a discretization $\hat{\pi}^\epsilon$ of $\hat{\pi}$ such that $\mathrm{d}(\hat{\pi}, \hat{\pi}^\epsilon  ) \leq \epsilon^{1/4}$. This approximation will introduce at most $L \epsilon^{1/4}$ error by the Lipschitz continuity of $f_M^1(\pi)$ and $f_M^2(\pi)$, so
\begin{equation*}
\liminf_{N \to \infty} F^\epsilon_N(\bvb, \bQ) \geq \liminf_{M \to \infty} \Big( f_M^1(\hat\pi^\epsilon) - f_M^2(\hat\pi^\epsilon) \Big) - L\delta  - L\epsilon^{1/4}.
\end{equation*}
These paths are now fixed, so we can now compute its limit as $M \to \infty$. Applying Lemma \ref{lem:sharpupbd} to decouple the constraint on $\bQ$ asymptotically shows
\begin{equation*}
\liminf_{M \to \infty} f_M^1(\hat\pi^\epsilon) \geq \inf_{\bL} \frac{1}{2} \bigg( \tr(\bL \bQ) - n - \log |\bL|  + (\bL^{-1}_0 \veh,\veh) +  \sum_{0 \leq k \leq r-1} \frac{1}{x_k}  \log \frac{|\bL_{k + 1}|}{|\bL_k|} \bigg),
\end{equation*}
where $(\bL_k)_{0 \leq k \leq r}$ are defined with respect to the sequences $(x_k)_{-1 \leq k \leq r}$ and $(\bQ_k)_{0 \leq k \leq r}$ encoded by $\hat\pi^\epsilon$.
By the recursive computations \eqref{eq:upbd_Y_recur},
\begin{equation*}
\lim_{M \to \infty} f_M^2(\hat\pi^\epsilon) = \sum_{0 \leq k \leq r-1} x_k \cdot \Sum \big( \btheta (\bQ_{k + 1}) - \btheta (\bQ_{k})\big).
\end{equation*}
Taking $\epsilon \to 0$ and consequently $\delta \to 0$ removes all the error terms, so we conclude
\begin{align}\label{eq:lowbound}
\lim_{\epsilon \to 0}\liminf_{N \to \infty} F^\epsilon_N(\bvb, \bQ) \geq \inf_{\bL,\pi} \sP_{\bvb, \bQ}(\pi, \bL).
\end{align}

\bibliographystyle{plain}

\end{document}